\documentclass[11pt]{amsart}
\makeatletter
\renewcommand{\section}{\@startsection{section}{1}{0pt}{20pt}{6pt}{\large\bfseries}}
\makeatother

\usepackage{enumerate}
\usepackage{amscd}

\numberwithin{equation}{section}

\theoremstyle{plain}
  \newtheorem{thm}{Theorem}[section]
  \newtheorem{prop}[thm]{Proposition}
  \newtheorem{lemma}[thm]{Lemma}
   \newtheorem{cor}[thm]{Corollary}
\theoremstyle{definition}
  \newtheorem{remark}[thm]{Remark}

\newcommand{\R}{\mathbb{R}}
\newcommand{\C}{\mathbb{C}}
\newcommand{\Q}{\mathbb{Q}}

\newcommand{\Qe}{\mathbf{E}}
\newcommand{\E}{\mathbb{E}}

\newcommand{\Az}{\mathbf{Q}}

\renewcommand{\P}{{\mathbb{P}}}

\newcommand{\Ip}{\mathcal{I}_{\kappa,\delta}}
\newcommand{\N}{\mathcal{N}}

\newcommand{\Ap}{\mathcal{K}_{\alpha,\gamma}}
\newcommand{\Id}[1]{{{\mathbb{I}}}_{\{#1\}}}
\newcommand{\fx}[2]{{}_#1\Psi_#2}
\renewenvironment{abstract}{%
  \ifx\maketitle\relax
  \fi
  \global\setbox\abstractbox=\vtop \bgroup
  \normalfont\Small
  \Small
  \list{}{\labelwidth
    \leftmargin \rightmargin\leftmargin
    \listparindent\normalparindent \itemindent
    \parsep
    
    }%
  \item[\noindent]%
}{%
  \endlist\egroup
  \ifx\@setabstract\relax \@setabstracta \fi}

\begin{document}
\bibliographystyle{plain}
\title{Law of the exponential functional of a new family of one-sided L\'evy processes via self-similar continuous state
branching processes with immigration and the $\fx{p}{q}$ Wright hypergeometric
functions}

\vspace{1cm}

\author{P. Patie$^{(1)}$}
\address{Department of Mathematical Statistics and Actuarial Science\\ University of Bern, Sidlerstrasse, 5, CH-3012 Bern, Switzerland}
\email{patie@stat.unibe.ch} \keywords{L\'evy processes,
self-similarity, continuous state branching processes with
immigration, Wright hypergeometric functions}
\begin{abstract}
We first introduce and derive some basic properties of a
two-parameters $(\alpha,\gamma)$ family of one-sided L\'evy processes, with
$1<\alpha<2$ and $\gamma>-\alpha$. Their Laplace exponents are given
in terms of the Pochhammer symbol as follows
\begin{eqnarray*}
\psi^{(\gamma)}(\lambda) &= &
c\left((\lambda+\gamma)_{\alpha}-(\gamma)_{\alpha}\right),
\quad \lambda \geq 0,
\end{eqnarray*}
where $c$ is a positive constant, $(\lambda)_\alpha=\frac{\Gamma(\lambda+\alpha)}{\Gamma(
\lambda)} $ stands for the Pochhammer symbol and $\Gamma$ for the
gamma function. These are a generalization of the Brownian
motion, since in the limit case $\alpha \rightarrow 2$, we end up to
Brownian motion with drift $\gamma +\frac{1}{2}$.
Then, we proceed by computing the density of the
law of the exponential functional associated to some elements of
this family (and their dual) and some transformations of these
elements. More specifically, we shall consider the L\'evy processes which admit the following Laplace exponent, for any $\delta>\frac{\alpha-1}{\alpha}$,
\begin{eqnarray*}
\psi^{(0,\delta)}(\lambda)&=&\psi^{(0)}(\lambda)-\frac{\alpha \delta}{\lambda+\alpha-1}\psi^{(0)}(\lambda),\quad   \lambda \geq 0.
\end{eqnarray*}
These densities are expressed in terms of the Wright hypergeometric function $\fx{1}{1}$. By means of probabilistic arguments, we derive some interesting properties enjoyed by this function. On the way we also characterize explicitly the semi-group of the family of self-similar positive Markov processes
associated, via the Lamperti mapping, to the family of L\'evy processes with Laplace exponent $\psi^{(0,\delta)}$.
\end{abstract}

\maketitle

\vspace{5mm}
\section{Introduction}
The exponential functional of L\'evy processes plays an important
role from both theoretical and applied perspectives. Indeed, it
appears in various fields such as  diffusion processes in random
environments, fragmentation and coalescence processes, the classical
moment problems, mathematical finance, astrophysics \ldots We refer
to the  paper of Bertoin and Yor \cite{Bertoin-Yor-05} for a
thorough survey on this topic and a description of cases when the
law of such functional is known explicitly. We also refer, in the
case of the Brownian motion with drift,  to the two survey papers of
Matsumoto and Yor \cite{Matsumoto-Yor-05-1} and
\cite{Matsumoto-Yor-05-2} where the law of the exponential
functional allows to characterize several interesting stochastic
processes and to develop stochastic analysis related to Brownian
motions on hyperbolic spaces. Finally, we mention that Bertoin et
al.~\cite{Bertoin-Biane-Yor-04} have expressed the law of the
exponential functional of a Poisson point process by means of
$q$-calculus and have derived an interesting connection to the
indeterminate moment problem associated to a log-normal random
variable.

In this paper, we start by introducing through its Laplace exponent and its characteristics triplet, a two $(\alpha,\gamma)$-parameters
family of spectrally negative L\'evy processes, with $1<\alpha\leq2$ and $\gamma>-\alpha$. The Laplace exponent has the following form
\begin{eqnarray} \label{eq:f1}
\psi^{(\gamma)}(\lambda) &= &
c\left((\lambda+\gamma)_{\alpha}-(\gamma)_{\alpha}\right), \quad
\lambda \geq 0,
\end{eqnarray}
where $c$ is a positive constant,
$(\lambda)_\alpha=\frac{\Gamma(\lambda+\alpha)}{\Gamma( \lambda)} $
stands for the Pochhammer symbol and $\Gamma$ for the Gamma
function. It turns out that this family possesses several
interesting properties. Indeed, their L\'evy measures behave around
$0$ as the L\'evy measure of a stable L\'evy process of index
$\alpha$. Moreover, a member of this family has its law at a fixed
time which belongs to the domain of attraction of a stable law of
index $\alpha$. They include the family of Brownian motion with
positive drifts  but also negative since they admit negative
exponential moments. After studying some further basic properties of
this family, we compute the density of the law of the exponential
functional associated to some elements of this family (and their
dual) and some transformations of these elements. More specifically,
we shall consider the L\'evy processes which admit the following
Laplace exponent, for any $\delta>\frac{\alpha-1}{\alpha}$,
\begin{eqnarray} \label{eq:f2}
\psi^{(0,\delta)}(\lambda)&=&\psi^{(0)}(\lambda)-\frac{\alpha
\delta}{\lambda+\alpha-1}\psi^{(0)}(\lambda),\quad   \lambda \geq 0.
\end{eqnarray}
The densities of the corresponding exponential functionals are
expressed in terms of the Wright hypergeometric function
$\fx{1}{1}$. As a specific instance, we obtain the inverse Gamma law
and hence recover Dufresne's result \cite{Dufresne-90} regarding the
law of the exponential functional of a Brownian motion with negative
drifts. As a limit case, we show that the family encompasses the
inverse Linnik law.

The path we follow to derive such a law is as follows. The first stage consists on characterizing
the family of L\'evy processes associated, via the Lamperti mapping, to the family of self-similar continuous
time branching processes with immigration (for short cbip).  We mention that this family includes the family of
squared Bessel processes and belong to the class of affine term structure models in mathematical finance,
 see the survey paper of Duffie et al.~\cite{Duffie-Filipovic-03}.
 At a second stage, by using well-know results on cbip, we derive the spatial Laplace transforms of the
  semi-groups of this later family. Then, by means of inversion techniques, we compute, in terms of  a power series,
  the density of these semi-groups. In particular, we  get an expression of their entrance laws in terms of the Wright
  hypergeometric function $\fx{1}{1}$. Finally, we end up our journey by related these entrance laws
  to the law of the exponential functionals associated to the family  of L\'evy processes \eqref{eq:f2}.

On the way, we get an expression, as a power series, for the density
of the semi-group of a family of sspMp. We mention, that beside the
family of Bessel processes, such a semi-group is only known
explicitly for the so-called  saw-tooth processes (a piecewise
linear sspMp) which were studied by Carmona et
al.~\cite{Carmona-Petit-Yor-01}. Therein, the authors use the
multiplicative kernel associated to a gamma random variable to
obtain an intertwining relationship between the semi-group of the
Bessel processes and the family of piecewise linear sspMp.

The remaining part of the paper is organized as follows. In the next
  Section, we gather some  preliminary results. In
particular, we introduce a transformation which leaves invariant the
class of sspMp without diffusion coefficients. We also provide the
detailed computation of an integral which will appear several times
in the paper. Section 3 concerns the definition and the study of
some basic properties of the new family of one-sided L\'evy
processes \eqref{eq:f1}.  Section 4 is devoted to the statement and
the proof of the law of the exponential functionals under study.
Section 5 contains  several remarks regarding some representations
of the special functions which appear in this paper. Finally Section
6 is a summarize of interesting properties enjoy by the Wright
hypergeometric functions.

\section{Recalls and preliminary results}

\subsection{Self-similar positive Markov processes}
Let $\zeta=(\zeta_t)_{t\geq0}$ denotes a real-valued L\'evy process
starting from $x$. Then, for any  $\alpha>0$, introduce the time
change process
\begin{eqnarray*}
X_t &=& e^{\zeta_{A_t}}
\end{eqnarray*}
where
\begin{eqnarray*}
A_t &=& \inf \{ s \geq 0; \Sigma_s := \displaystyle{\int_0^s}
e^{\alpha \zeta_u} \: du
> t \}.
\end{eqnarray*}
Lamperti \cite{Lamperti-72} showed that the process $(X_t)_{t\geq0}$ is a
positive-valued $\frac{1}{\alpha}$ semi-stable Markov process (for
short sspMp) starting from $e^x$. That is, if, for $x>0$, $\Q_x$
denotes the law of  $X$ starting from $x$, then $X$ is a positive
valued Markov process which enjoys the $\alpha$-self-similarity
property, i.e.~for any $c>0$, we have the equality in distribution
\begin{equation}
\left((X_{c^{\alpha}t})_{t\geq0}, \Q_{cx}\right)
\stackrel{(d)}{=}\left((cX_{t})_{t\geq0}, \Q_x\right).
\end{equation}
Lamperti showed that the mapping above is actually one to one. We
shall refer to this time change transformation as the Lamperti
mapping. We shall use the symbol $\Qe$ for the expectation operator
associated to $\Q$. Finally, let us denote by $T_0^X=\inf\{s\geq0; \: X_{s-}=0,X_s=0\}$ the lifetime of $X$.
Note that $T_0^X=\infty$ or $<\infty$ a.s.~according to $E[\zeta_1]\geq0$ or not. We proceed by stating a  general result on a
transformation between sspMp.
\begin{prop}
Let $\zeta$, with characteristic exponent $\Upsilon$, be the image
via the Lamperti mapping of an $\alpha$-sspMp $X$ and fix $\beta>0$.
Then, if $\sigma=0$, the $\frac{\alpha}{\beta}$-sspMp $X^{\beta}$
and the $\alpha$-sspMp $X^{(\beta)}$, obtained as the image via the
Lamperti mapping of the L\'evy process $\beta \zeta$, have the same
image via this mapping. Otherwise, the characteristics exponent of
 the image
via the Lamperti mapping of the $\frac{\alpha}{\beta}$-sspMp $X^{\beta}$ is
\begin{eqnarray*}
 \Upsilon(\beta \lambda) +
\frac{\sigma}{2}\beta(\beta-1)\lambda,\quad \lambda \in \R.
\end{eqnarray*}
\end{prop}
The proof is split into the following two lemmas.
\begin{lemma} \label{lem:ts}
For any $\beta>0$, we have the following relationship between sspMp,
\begin{equation} \label{eq:rss}
X^{(\beta)}_t = X^{\beta}_{\int_0^t
X^{(\beta)\alpha(\frac{1}{\beta}-1)}_sds}, \: t<T_0^X,
\end{equation}
where $(X^{(\beta)},\Q)$ is the sspMp associated to the L\'evy
process $\beta \zeta$.
\end{lemma}
\begin{proof}
Since $\beta>0$, note that the lifetime of both processes are either
finite a.s or infinite a.s. Then, set $A_t=\int_0^tX^{-\alpha}_sds$
and observe that the Lamperti mapping yields, for any $t<T_0^X$,
\begin{eqnarray*}
\log\left(X^{\beta}_t\right) &=&
\beta\zeta_{A_t}\\  &=& \log\left(X^{(\beta)}_{\int_0^{A_t}e^{\alpha \beta \xi_s}}ds\right) \\
&=& \log\left(X^{(\beta)}_{\int_0^{t}X^{\alpha(\beta-1)}_sds}\right)
\end{eqnarray*}
which completes the proof.
\end{proof}
Next, we recall that the infinitesimal generator, $\mathbf{Q}$, of
the $\alpha$-sspMp $X$ is given for any $\alpha>0$ and $f$ a smooth
function on $\R^+$, by
\begin{eqnarray*}
\mathbf{Q} f(x)&=& x^{-\alpha} \left(\frac{\sigma}{2} x^{2}f''(x) +
b x f'(x)\right.
\\&+&
\left.\int^{+\infty}_{-\infty}\left((f(e^{r}x)-f(x))-xf'(x)r\Id{r\leq1}\right)\nu(dr)\right)
\nonumber
\end{eqnarray*}
where the three parameters $ b \in \R, \: \sigma \geq 0$ and the
measure $\nu$ which satisfies the integrability condition
$\int_{-\infty}^{\infty} (1\wedge |r|^2 )\: \nu(dr) <+ \infty$ form
the characteristic triplet of the L\'evy process $\zeta$, the image
of $X$ by the Lamperti mapping.
\begin{lemma} \label{lem:ps}
With the notation as above, for any $\beta>0$, the laplace exponent of the L\'evy process
associated to the $\frac{\alpha}{\beta}$-sspMp $\left(X^{\beta},
\Q\right)$ is given by
\begin{eqnarray*}
\Upsilon(\beta \lambda ) +
\frac{\sigma}{2}\beta(\beta-1)\lambda,\quad \lambda \in \R.
\end{eqnarray*}
\end{lemma}
\begin{proof}
Let us denote by $\Az^{\beta}$ the infinitesimal generator of
$(X^{\beta},\Q)$ and set $p_{\beta}(x)=x^{\beta}$ for any $\beta>0$.
Since the function $x\mapsto p_{\beta}(x)$ is a homeomorphism of
$\R^+$ into itself, we obtain, after some easy computations,
\begin{eqnarray*}
\Az^{\beta} f(x) &= &\Az(f\circ p_{\beta})(x^{1/\beta}) \\
& = &   x^{-\alpha/\beta} \left(\frac{\sigma}{2} \beta^2x^{2}f''(x)
+  \left(b\beta+ \frac{\sigma}{2}\beta (\beta-1)\right)x
f'(x)\right.
\\&+&
\left.\int^{+\infty}_{-\infty}\left((f(e^{\beta r}x)-f(x))-\beta
xf'(x)r\Id{r\leq1}\right)\nu(dr)\right).
\end{eqnarray*}
The proof is completed by identification.
\end{proof}
The proof of the Proposition follows by putting pieces together.

\subsection{An integral computation}
We proceed by computing an integral, expressed in terms of ratios of
gamma functions, which will be very important for the sequel.
Indeed, it will be useful for computing the characteristic triplet
associated to the families \eqref{eq:f1} and \eqref{eq:f2}, which
will be provided in Section \ref{sec:l}.  We fix the following
constants
\begin{eqnarray*}
c&=& -\frac{1}{\alpha \cos\left(\frac{\alpha
\pi}{2}\right)} >0\\
c_{\alpha}&=&\frac{c}{\Gamma(-\alpha)} >0.
\end{eqnarray*}
The motivation for the choice of theses constants is given in Remark
\ref{rem:l} below.
\begin{thm} For any $\alpha,\lambda,\gamma \in \C$, with $1<\Re(\alpha)<2$, $\Re(\lambda)>0$ and $\Re(\alpha+\gamma)> 0$, we have
\begin{equation} \label{eq:int_poch}
c_{\alpha}\int_0^{1}\frac{(u^{\lambda}-1)u^{\alpha+\gamma-1}-\lambda(u-1)}{(1-u)^{\alpha+1}}du
=
c\left((\lambda+\gamma)_\alpha-(\gamma)_\alpha\right)-\frac{c_{\alpha}\lambda}{\alpha-1},
\end{equation}
where we recall that
$(\lambda)_{\alpha}=\frac{\Gamma(\lambda+\alpha)}{\Gamma(\lambda)}$
is the Pochhammer symbol and $c_{\alpha}=\frac{c}{\Gamma(-\alpha)}
>0$.
\end{thm}
\begin{proof}
In the sequel, we denote by $\mathcal{F}_{\lambda,\alpha,\gamma}$
the integral of the left-hand side on \eqref{eq:int_poch} divided by $c_{\alpha}$. Before
starting the proof, we recall the following integral representation
of the Beta function, see e.g.~Lebedev \cite{Lebedev-72},
\begin{eqnarray*}
\mathcal{B}(\gamma,\alpha)&=& \frac{\Gamma(\gamma)\Gamma(\alpha)}{\Gamma(\gamma+\alpha)}\\
&=&\int_0^{1}(1-v)^{\gamma-1}v^{\alpha-1}dv, \quad \Re(\gamma)>0,\:
\Re(\alpha)>0,
\end{eqnarray*}
and the recurrence relation for the gamma function
$\Gamma(\lambda+1)=\lambda\Gamma(\lambda)$. Then, reiteration of
integrations by parts yield
\begin{eqnarray*}
\mathcal{F}_{\lambda,\alpha,\gamma} &=&
-\frac{\lambda(\lambda+\gamma)_{\alpha}\Gamma(1-\alpha)}{\alpha(\lambda+\gamma+\alpha-1)}+\frac{\lambda}{\alpha-1}
-\frac{\alpha+\gamma-1}{\alpha}\mathcal{F}
\end{eqnarray*}
where we have set
\begin{eqnarray*}
\mathcal{F}&:=&\int_0^{1}(u^{\lambda}-1)u^{\alpha+\gamma-2}(1-u)^{-\alpha}du
\end{eqnarray*}
and we have used the condition $\Re(\alpha+\gamma)> 0$.
Next, according to the binomial expansion, we have
\begin{eqnarray*}
\mathcal{F}&=&\sum_{n=0}^{\infty}\frac{(\alpha)_{n}}{n!}\int_0^{1}(u^{\lambda}-1)u^{n+\alpha+\gamma-2}du\\
&=&-\frac{1}{\lambda}\sum_{n=0}^{\infty}\frac{(\alpha)_{n}}{n!}\int_0^{1}(1-v)v^{\frac{n+\alpha+\gamma-1}{\lambda}-1}dv\\
&=&-\lambda
\sum_{n=0}^{\infty}\frac{(\alpha)_{n}}{(n+\alpha+\gamma-1+\lambda)(n+\alpha+\gamma-1)n!}\\
&=&-\frac{\lambda}{\alpha+\gamma-1}
\sum_{n=0}^{\infty}\frac{(\alpha)_{n}(\alpha+\gamma-1)_{n}}{(\alpha+\gamma)_{n}n!(n+\alpha+\gamma-1+\lambda)}.
\end{eqnarray*}
Using the power series of the ${}_2\mathcal{F}_{1}$ hypergeometric
functions, see e.g.~\cite[9.1]{Lebedev-72},
\begin{eqnarray*}
{}_2\mathcal{F}_{1}(\alpha,\beta;\gamma;z)&=&\sum_{n=0}^{\infty}\frac{(\alpha)_{n}(\beta)_n}{(\gamma)_n n!}z^n, \quad \mid z\mid <1,\\
\end{eqnarray*}
we observe that
\begin{eqnarray*}
\mathcal{F}&=&-\frac{\lambda}{\alpha+\gamma-1}\lim_{z\rightarrow 1^-}\int_0^z u^{\lambda+\alpha+\gamma-2} {}_2\mathcal{F}_{1}(\alpha,\alpha+\gamma-1;\alpha+\gamma;u)du\\
&=&-\frac{\lambda}{\alpha+\gamma-1}\frac{\Gamma(\lambda+\alpha+\gamma-1)}{\Gamma(\lambda+\alpha+\gamma)}\lim_{z\rightarrow
1^-}{}_3\mathcal{F}_{2}(^{\alpha,\alpha+\gamma-1,\lambda+\alpha+\gamma-1}_{\alpha+\gamma,\lambda+\alpha+\gamma};z)
\end{eqnarray*}
where the last line follows from
\cite[7.512(5)]{Gradshteyn-Ryzhik-00} and ${}_3\mathcal{F}_{2}$ is the hypergeometric function of degree $(3,2)$. Note that this later
representation holds for $\Re(\alpha)<2$. We proceed by using a
 result of Milgram \cite[(11)]{Milgram-06} regarding the
limit of the ${}_3\mathcal{F}_2$ function which is
\begin{eqnarray*}
\lim_{z\rightarrow
1^-}{}_3\mathcal{F}_{2}(^{\alpha,\alpha+\gamma-1,\lambda+\alpha+\gamma-1}_{\alpha+\gamma,\lambda+\alpha+\gamma};z)&=&\frac{-(\alpha+\gamma-1)\Gamma(\alpha+\gamma+\lambda)\Gamma(1-\alpha)}{\lambda\Gamma(\gamma+\lambda)}\\
&+&\frac{(\alpha+\gamma+\lambda-1)\Gamma(\alpha+\gamma)\Gamma(1-\alpha)}{\Gamma(\gamma)\lambda}.
\end{eqnarray*}
It follows that
\begin{eqnarray} \label{eq:if}
\mathcal{F}&=&\frac{\lambda\Gamma(1-\alpha)}{(\gamma+\alpha-1)(\lambda+\gamma+\alpha-1)}\left(\frac{(\gamma+\alpha-1)(\lambda+\gamma)_{\alpha}}{\lambda}-\frac{(\lambda+\gamma+\alpha-1)(\gamma)_\alpha}{\lambda}\right)\nonumber\\
&=&\Gamma(1-\alpha)\left(\frac{(\lambda+\gamma)_{\alpha}}{\lambda+\gamma+\alpha-1}-\frac{(\gamma)_\alpha}{(\gamma+\alpha-1)}\right).
\end{eqnarray}
Finally, we obtain
\begin{eqnarray*}
\mathcal{F}_{\lambda,\alpha,\gamma}&=&\frac{\Gamma(1-\alpha)}{\alpha}\left(-\frac{(\lambda+\gamma)_{\alpha}(\lambda+\gamma+\alpha-1)}{\lambda+\gamma+\alpha-1}
+(\gamma)_\alpha \right)+\frac{\lambda}{\alpha-1}
\end{eqnarray*}
which completes the proof by recalling that $c_{\alpha}=-c\frac{\Gamma(1-\alpha)}{\alpha}$.
\end{proof}

\section{Basic properties of the family
$\left(\xi,\P^{(\gamma)}\right)$} \label{sec:l} Let us denote by
$\P^{(\gamma)}=\left(\P^{(\gamma)}_x\right)_{x\in \R}$ the family of
probability measures of the process $\xi$ such that
$\P^{(\gamma)}_x(\xi_0=x)=1$ and  recall that the Laplace exponent,
denoted by $\psi^{(\gamma)}$, of the process $(\xi,\P^{(\gamma)})$
has the following form
\begin{eqnarray*}
\psi^{(\gamma)}(\lambda) &= &
c\left((\lambda+\gamma)_{\alpha}-(\gamma)_{\alpha}\right), \quad
\lambda \geq 0,
\end{eqnarray*}
where  $(\lambda)_\alpha=\frac{\Gamma(\lambda+\alpha)}{\Gamma(
\lambda)} $ stands for the Pochhammer symbol and  $c= -\frac{1}{\alpha \cos\left(\frac{\alpha
\pi}{2}\right)}$, and the parameters $\alpha,\gamma$
belong to the set
\begin{eqnarray*}
\Ap= \left\{1 < \alpha < 2, \gamma >-\alpha\right\}.
\end{eqnarray*}
We denote by $\E^{(\gamma)}_x$ the expectation operator associated
to $\P^{(\gamma)}_x$ and write simply $\E^{(\gamma)}$ for
$\P^{(\gamma)}_0$. In what follow, we show that it is the Laplace
exponent of a spectrally negative L\'evy process, we also provide
its characteristic triplet and derive some basic properties.
\begin{prop} \label{prop:nl}
\begin{enumerate}
\item For $\alpha,\gamma \in \Ap$, the process $\left(\xi,
\P^{(\gamma)}\right)$ is a spectrally negative L\'evy process with
finite quadratic variation.  More specifically, we have
\begin{eqnarray*}
\psi^{(\gamma)}(\lambda) &= & \tilde{c}_{\alpha} \lambda +\int_{-\infty}^0\left(e^{\lambda y}-1-\lambda y\Id{\mid y\mid<1}\right)\nu(dy)
\end{eqnarray*}
where
\begin{equation*}
\nu(dy)=c_{\alpha}\frac{e^{(\alpha+\gamma)y}}{(1-e^{y})^{\alpha+1}}dy,
\quad y<0,
\end{equation*}
and
\begin{eqnarray*}
 \tilde{c}_{\alpha}&=&
c_{\alpha}\sum_{k=1}^{\infty}\frac{1}{k(k-\alpha)}\left(\mathcal{B}\left(\frac{e-1}{e};k+1-\alpha,\alpha+\gamma-1\right) -\right.\\
& &\left.\left(\frac{e-1}{e}\right)^{k-\alpha}((\alpha+\gamma-1)e^{1-\alpha-\gamma}-1)\right)
\end{eqnarray*}
where $\mathcal{B}\left(.;.,.\right)$ stands for the incomplete Beta function.
\item  The random variable $\left(\xi_1,
\P^{(\gamma)}\right)$ admits negative exponential moments of order
lower than $\gamma+\alpha$, i.e. for any $\lambda<\gamma+\alpha$, we
have
\begin{eqnarray*}
\E^{(\gamma)}\left[e^{-\lambda \xi_1} \right]<+\infty.
\end{eqnarray*}
\item We have  the following invariance property by Girsanov transform,
i.e.~ for any $\gamma \in \Ap$,
\begin{equation} \label{eq:esscher}
d\P_0^{(\gamma)}{}_{\mid F_t}=e^{\gamma
\xi_t-\psi^{(0)}(\gamma)t}d\P_0^{(0)}{}_{\mid F_t}, \quad t>0.
\end{equation}
\item The first moments have the following expressions
\begin{eqnarray} \label{eq:first_moment}
\E^{(\gamma)}\left[\xi_1 \right] &=& c
(\gamma)_{\alpha}(\Psi(\gamma+\alpha)-\Psi(\gamma)) \:, \gamma>-\alpha ,\nonumber \\
\E^{(0)}\left[\xi_1 \right] &=& c\Gamma(\alpha)\\
\E^{(1-\alpha)}\left[\xi_1 \right] &=& c\frac{1}{\Gamma(1-\alpha)}\left(-E_{\gamma}-\Psi(1-\alpha)\right)\\
\E^{(-1)}\left[\xi_1 \right] &=&- c \Gamma(\alpha-1)
\end{eqnarray}
where $\Psi(\lambda)=\frac{\Gamma'(\lambda)}{\Gamma(\lambda)}$ is
the digamma function and $E_{\gamma}$ stands for Euler-Mascheroni
constant.
\item For any
$1<\alpha <2$, there exits $\gamma_{\alpha}$, with $-\alpha
<\gamma_{\alpha} <0$, such that $\E^{(\gamma_{\alpha})}\left[\xi_1
\right]=0$. Moreover, for any $\gamma<\gamma_{\alpha}$ the Cram\'er
condition holds, i.e.~for any $\gamma<\gamma_{\alpha}$, there
exists $\lambda_{\alpha}>0$ such that
\begin{eqnarray}
\E^{(\gamma)}\left[e^{\lambda_{{\alpha}} \xi_1} \right] &=&1.
\end{eqnarray}
\item Consequently, for any fixed $1<\alpha<2$, the process $\left(\xi_t,
\P^{(\gamma_{\alpha})}\right)$ oscillates and otherwise
$\lim_{t\rightarrow \infty} \left(\xi_t, \P^{(\gamma)}\right) =
\textrm{sgn}(\gamma-\gamma_{\alpha})\infty$ a.s.
\item For $\gamma=0$ or $-1$, the
scale function of  $\left(\xi, \P^{(\gamma)}\right)$ is given by
\begin{equation*}
\mathcal{W}^{(\gamma)}(x)=\frac{1}{ \Gamma(\alpha)}e^{-\gamma
x}(1-e^{-x})^{\alpha-1},\quad x>0.
\end{equation*}
\item Finally, we have the following two limits results:
\begin{enumerate}
\item \label{it:b} For any $\gamma \in \Ap$, the process $\left(\xi,
\P^{(\gamma)}\right)$  converges in distribution as $\alpha
\rightarrow 2$ to a  Brownian motion with drift
$\gamma+\frac{1}{ 2}$.
\item \label{it:s}For any fixed $1<\alpha<2$, the process $\left(\lambda^{\frac{1}{\alpha}}\xi_{\lambda t},
\P^{(0)}\right)$  converges in distribution as $\lambda \rightarrow \infty$ to  a spectrally negative $\alpha$-stable process.
\end{enumerate}
\end{enumerate}
\end{prop}
\begin{proof}
\begin{enumerate}
\item
First, from \eqref{eq:int_poch}, we deduce that
\begin{eqnarray*}
\psi^{(\gamma)}(\lambda) &= &
c\left((\lambda+\gamma)_{\alpha}-(\gamma)_{\alpha}\right) \\
&= &  \frac{c_{\alpha}\lambda}{1-\alpha} + \int_{0}^1
\left((u^{ \lambda}-1)u^{\alpha+\gamma-1}- \lambda (u-1)\right)\frac{c_{\alpha}du}{(1-u)^{\alpha+1}}\\
&= & \left(c_{\alpha}\left(\int_{0}^{1} \frac{(\log(u)-(u-1))\Id{
\mid
 \log(u)\mid
<1}}{(1-u)^{\alpha+1}}du+\int^{\frac{1}{e}}_{-\infty}(1-u)^{-\alpha}du\right)
\right.\\
&+& \left. c_{\alpha}\int_{0}^1
\frac{\left(u^{\alpha+\gamma-1}-1\right)\log(u)\Id{ \mid
 \log(u)\mid
<1}}{(1-u)^{\alpha +1}}du \right)\lambda
\\&+& \int_{0}^1
\left(u^{ \lambda}-1-\lambda \log(u) \Id{\mid \log(u)\mid
<1}\right)\frac{c_{\alpha}u^{\alpha+\gamma-1}du}{(1-u)^{\alpha+1}}\\
&= & \tilde{c}_{\alpha}\lambda + \int_{0}^{\infty} \left(e^{\lambda
y}-1-\lambda y \Id{\mid y \mid
<1}\right)\frac{c_{\alpha}e^{(\alpha+\gamma)y}du}{(1-e^{y})^{\alpha+1}}.
\end{eqnarray*}
where we have set
\begin{eqnarray*}
 \tilde{c}_{\alpha}=c_{\alpha}\left(\int_{0}^{1} \frac{(\log(u)-(u-1))\Id{
\mid
 \log(u)\mid
<1}}{(1-u)^{\alpha+1}}du+\int^{\frac{1}{e}}_{-\infty}(1-u)^{-\alpha}du \right.\\
\left. + \int_{0}^1
\frac{\left(u^{\alpha+\gamma-1}-1\right)\log(u)\Id{ \mid
 \log(u)\mid
<1}}{(1-u)^{\alpha +1}}du \right)
\end{eqnarray*}
Hence, we recognize the L\'evy-Khintchine representation of a
one-sided L\'evy process. Moreover, the quadratic finite variation
property follows from the following asymptotic behavior of the Pochhammer symbol, see e.g.~\cite{Lebedev-72},
\begin{eqnarray} \label{eq:ag}
(z+ \gamma)_{\alpha}^{-1} &=&
z^{-\alpha}\left[1+\frac{(-\alpha)(2\gamma-\alpha-1)}{2z} +
O(z^{-2})\right],\: \mid \textrm{arg } z\mid < \pi-\delta,\: \delta>0,
\end{eqnarray}
and the condition $1<\alpha<2$.
It remains to compute the constant
$\tilde{c}_{\alpha}$. Performing the change of variable $v=1-u$, we get, for the first term on the left-hand side of the previous identity,
\begin{eqnarray*}
\int^{\frac{e-1}{e}}_0 (\log(1-v)+v)\frac{c_{\alpha}}{ v^{\alpha +1}}dv   &=& -
c_{\alpha}\sum_{k=2}^{\infty}\frac{1}{k}\int^{\frac{e-1}{e}}_0
 v^{k-\alpha-1}dv \\
&=&-
c_{\alpha}\sum_{k=2}^{\infty}\frac{1}{k(k-\alpha)}\left(\frac{e-1}{e}\right)^{k-\alpha}.
\end{eqnarray*}
Moreover, proceeding as above, we have
\begin{eqnarray*}
 & & \int_{0}^1
\frac{\left(u^{\alpha+\gamma-1}-1\right)\log(u)\Id{ \mid
 \log(u)\mid
<1}}{(1-u)^{\alpha +1}}du \\&=& -
c_{\alpha}\sum_{k=1}^{\infty}\frac{1}{k}\int^{\frac{e-1}{e}}_0
 v^{k-\alpha-1}((1-v)^{\alpha+\gamma-1}-1)dv \\
&=&
c_{\alpha}\sum_{k=1}^{\infty}\frac{\alpha+\gamma-1}{k(k-\alpha)}\left(\int^{\frac{e-1}{e}}_0
 v^{k-\alpha}(1-v)^{\alpha+\gamma-2}dv \right.\\
& & \left. -\left(\frac{e-1}{e}\right)^{k-\alpha}e^{1-\alpha-\gamma}\right)\\
&=&
c_{\alpha}\sum_{k=1}^{\infty}\frac{\alpha+\gamma-1}{k(k-\alpha)}\left(\mathcal{B}\left(\frac{e-1}{e};k+1-\alpha,\alpha+\gamma-1\right) -\right.\\
& &\left.\left(\frac{e-1}{e}\right)^{k-\alpha}e^{1-\alpha-\gamma}\right).
\end{eqnarray*}

Putting pieces together, one gets
\begin{eqnarray*}
 \tilde{c}_{\alpha}&=&
c_{\alpha}\sum_{k=1}^{\infty}\frac{1}{k(k-\alpha)}\left(\mathcal{B}\left(\frac{e-1}{e};k+1-\alpha,\alpha+\gamma-1\right) -\right.\\
& &\left.\left(\frac{e-1}{e}\right)^{k-\alpha}((\alpha+\gamma-1)e^{1-\alpha-\gamma}-1)\right)
\end{eqnarray*}
which completes the description of the characteristic triplet of
$\left(\xi, \P^{(\gamma)}\right)$.
\item This item follows from the fact that the mapping $\lambda\rightarrow \psi^{(\gamma)}(\lambda)$ is well defined on $(-\gamma+\alpha,\infty)$.
\item It is simply the Esscher transform.
\item The expressions
for the first moment
 of $\left(\xi_1, \P^{(\gamma)}\right)$
is obtained from the formula
\begin{eqnarray}
\frac{\partial}{\partial{\lambda}}(\lambda+\gamma)_{\alpha}=(\lambda+\gamma)_{\alpha}(\Psi(\lambda+\gamma+\alpha)-\Psi(\lambda+\gamma)).
\end{eqnarray}
 Moreover, for $\gamma=0$ and $\gamma=-1$
we use the recurrence formulae for the digamma function, $\Psi(u+1)=\frac{1}{u}+\Psi(u)$, see
\cite[Formula 1.3.3]{Lebedev-72} and for the Gamma function.
\item
 Since, for any $\gamma \in \Ap$, the mapping $\lambda \mapsto \psi^{(\gamma)}(\lambda)$ is convex and continuously differentiable on $(\alpha-\gamma,\infty)$, its derivative is continuous and increasing.
 Hence, the mapping $\gamma  \mapsto \psi^{(\gamma)}(0^+)$ is continuous and increasing on $\Ap$. Moreover,  noting that
 $\E^{(-1)}\left[\xi_1 \right]<0<\E^{(0)}\left[\xi_1 \right]$ for any $1<\alpha<2$,
 we deduce that for each $1<\alpha<2$, there exists an unique $-1<\gamma_{\alpha}<0$ such that
 $\E^{(\gamma_{\alpha})}\left[\xi_1 \right]=0$. Since for any $\gamma<\gamma_{\alpha}$, $ \E^{(\gamma)}\left[\xi_1 \right]$
 is negative and $\lambda  \mapsto \psi^{(\gamma)}(\lambda)$ is convex with $\lim_{\lambda \rightarrow \infty} \psi^{(\gamma)}(\lambda) = \infty$,
 we deduce that there exists $\lambda_{\alpha}>0$ such that $\psi^{(\gamma)}(\lambda_{\alpha})=0$.

 \item The long time behavior of  $\left(\xi, \P^{(\gamma)}\right)$ follows from the previous item and the strong law of large numbers.
\item The expression of the scale function
is derived from the following identity, see
\cite[3.312,1.]{Gradshteyn-Ryzhik-00}, for
$\Re(\alpha),\Re(\lambda+\gamma)>0, $
\begin{equation*}
\int_0^{\infty} e^{-(\lambda+\gamma) x}(1-e^{-x})^{\alpha-1}dx =
\frac{\Gamma(\lambda+\gamma)\Gamma(\alpha)}{\Gamma(\lambda+\gamma
+\alpha)}.
\end{equation*}
\item It is enough to show that the random variable $\left(\xi_1,
\P^{(\gamma)}\right)$ converges in law to a normal distribution. By
continuity of the function $\psi^{(\gamma)}(\lambda)$ in $\alpha$,
we get the result after easy manipulation of the Gamma function.
Similarly, for the second limit, we simply need to show that the
random variable $\left( \xi_1, \P^{(0)}\right)$ belong to the domain
of attraction of  a stable distribution, i.e. $\lim_{\eta
\rightarrow \infty} \eta \psi^{(0)}(\eta^{-\frac{1}{\alpha}}\lambda)
= c\lambda^{\alpha}$. The claim follows by means of the asymptotic of
the ratio of Gamma functions, see \eqref{eq:ag}.
\end{enumerate}
\end{proof}
\begin{remark}
 The choice of the constant $c_{\alpha}$ is motivated by the item \eqref{it:s}. Actually, the coefficients of
 the Laplace exponent of a completely asymmetric stable random variable is given for any $\iota>0$ by $c(\iota)=
-\frac{\iota^{\alpha}}{\alpha \cos\left(\frac{\alpha \pi}{2}\right)}
>0$, see e.g.~\cite[Proposition 1.2.12]{Sam-Taq-94}. Thus, we have made the choice $\iota^{\alpha}=1$.
\end{remark}
 \begin{remark} \label{rem:l} Note that Caballero and Chaumont \cite{Caballero-Chaumont-06} show, by characterizing its characteristic triplet, that the process $(\xi,\P^{(0)})$ is the image via the Lamperti mapping of a spectrally negative regular $\alpha$-stable process conditioned to stay positive. For any $\gamma \in \Ap$, the expression of the Laplace exponent appears as an example but without proof, in \cite{Patie-06c}.
\end{remark}

\section{Law of some exponential functionals via$\ldots$}
In this section, we derive the explicit law of the exponential
functionals associated to some elements of the family of L\'evy
processes introduced in this paper and to some transforms of these
elements.

We proceed by defining a family of functions which will appear
several times in the remaining part of the paper. The Wright
hypergeometric function is defined as, see \cite[1,Section
4.1]{Erdelyi-55},
\begin{equation*}
 \fx{p}{q}\left( \left.\begin{array}{c}
                  (A_1,a_1),\ldots(A_p,a_p) \\
                  (B_1,b
                  _1)\ldots(B_q,b_q)
                \end{array} \right|
 \: z\right)
= \sum_{n=0}^{\infty} \frac{\prod_{i=1}^p\Gamma(A_i n+a_i
)}{\prod_{i=1}^q\Gamma(B_in +b_i)}\frac{z^n}{n!}
\end{equation*}
where $p,q$ are nonnegative integers, $a_i \in \C \: (i=1\ldots
p), b_j\in \C \: (j=1\ldots q)$ and the coefficients $A_i \in
\R^+\: (i=1\ldots p)$ and $B_j \in \R^+ \:(j=1\ldots q)$ are such
that
$1+\sum_{i=1}^q B_i-\sum_{i=1}^p A_i \geq 0.$
Under such conditions, it follows from the following asymptotic
formula of ratio of gamma functions, see \eqref{eq:ag},
\begin{eqnarray*}
\frac{\Gamma(z+ \gamma)}{\Gamma(z+\gamma + \alpha)} &=&
z^{-\alpha}\left[1+\frac{(-\alpha)(2\gamma-\alpha-1)}{2z} +
O(z^{-2})\right], \quad  \mid \textrm{arg } z\mid < \pi-\delta,
\end{eqnarray*}
that  $\fx{p}{q}(z)$ is an entire function with respect to $z$. We
postpone to Section \ref{sec:w} for a more detailed description of
this class of function.
\subsection{$\ldots$ an important continuous state branching process}
In this part, we establish a connection between the dual of an
element of the family of L\'evy processes introduced in this paper
and a self-similar continuous state branching process. As a
byproduct, we derive the explicit law of the exponential functional
associated to the dual of this element.

To this end,  we recall that continuous
state branching processes form a class of non negative valued Markov
processes which appear as limit of integer valued branching
processes.  Lamperti \cite{Lamperti-67b} showed that when the units
and the initial state size is
 allowed to tend to infinity, the limiting process is necessarily a
self-similar continuous state branching process with index lower than $1$. We denote  this process by
$(Y,\Q^{(0)})$, i.e.~$\Q^{(0)}=(\Q^{(0)})_{x>0}$ is a family of probability measures such that $\Q^{(0)}_x(Y_0=x)=1$. The associated expectation operator is $\Qe^{(0)}$. Moreover, Lamperti showed
that the semi-group of $(Y,\Q^{(0)})$ is characterized by its spatial Laplace
transform as follows, for any $ \lambda,x \geq 0$,
\begin{eqnarray} \label{eq:lbs}
\Qe^{(0)}_x\left[e^{-\lambda Y_t}\right]= e^{-xd\lambda \left(1+ct
\lambda^{\kappa}\right)^{-1/\kappa}},\
\end{eqnarray}
for some $0<\kappa\leq 1$ and some positive constants $d,c$. On the
other hand, he also observed in \cite{Lamperti-67a} that a
continuous state branching process can be obtained from a \emph{spectrally
positive} L\'evy process, $\zeta$, by a time change. More precisely,
consider $\zeta$ started at $x>0$ and write
\begin{eqnarray*}
T^\zeta_0&=&\inf\{s\geq 0;\: \zeta_s=0\}.
\end{eqnarray*}
Next, let
\begin{eqnarray*}
\Lambda_t &=& \int_0^t \zeta^{-1}_s ds,  \quad t<T^\zeta_0
\end{eqnarray*}
and
\begin{eqnarray*}
  V_t &=& \inf\{s\geq
0; \: \Lambda_s>t\} \wedge T^\zeta_0.
\end{eqnarray*}
Then,  the time change process $\zeta\circ V$ is a continuous state
branching process starting at $x$ and the Laplace exponent $\varphi$
of $\zeta$ is called the branching mechanism.  Finally, we recall
that the law of the  absorption time
\begin{eqnarray*}
i&=&\inf\{s\geq 0; \zeta\circ V_s=0\}
\end{eqnarray*}
has been computed explicitly by Grey \cite{Grey-74}. More
specifically, put $\phi(0)=\inf\{s\geq0;\varphi(s)>0\}$ (with the
usual convention that $\inf \emptyset =\infty$) and  assume that
$\int^{\infty}\varphi^{-1}(s)ds<\infty$ and $\phi(0)<\infty$, then
\begin{eqnarray} \label{eq:li}
\textrm{the law of } g(i) \textrm{ for }{\zeta\circ V} \textrm{ starting at } x>0 \textrm{ is   exponential with parameter } x,
\end{eqnarray}
where $g:(0,\infty)\rightarrow (\phi(0),\infty)$ is the inverse
mapping of $\int_t^{\infty}\varphi^{-1}(s)ds$.
 We are now ready to
state  and proof the main result of this part.
\begin{thm} \label{thm:e1}
For any $0<\kappa\leq1$, we have
\begin{eqnarray*}
\left(\int_0^{\infty}e^{- \kappa\xi_s }ds, \P_0^{(0)}\right) &
\sim &
c G(1)^{-\kappa}
\end{eqnarray*}
where $G(1)$ is a Gamma random variable of parameter $1$.
\end{thm}
In the case $\kappa=1$,  $\left(\xi, \P^{(0)}\right)$ is a Brownian motion with positive drift $\frac{1}{2}$, see item \ref{it:b} in Proposition \ref{prop:nl}, and the result above corresponds to Dufresne's result \cite{Dufresne-90}.
We split the proof in several lemmas. First, let us denote by $X$  a
\emph{spectrally positive} $\alpha$-stable
L\'evy process killed when it hits zero. It is plain that it is also
a $\alpha$-sspMp. We have the following.

\begin{lemma} \label{lem:fd}
The L\'evy process associated, via the Lamperti mapping, to
$X$ is the dual, with respect to the Lebesgue measure,
of $(\xi,\P^{(0)})$, i.e.~$(-\xi,\P^{(0)})$.
\end{lemma}
\begin{proof}
First, by Hunt switching identity, see e.g.~Getoor and Sharpe \cite{Getoor-Sharpe-84}, we have that $X$
is in duality, with respect to
the Lebesgue measure, with the spectrally negative $\alpha$-stable L\'evy process killed when entering the negative
real line. Moreover, it is well known, see e.g.~Bertoin \cite{Bertoin-96} that the spectrally negative $\alpha$-stable
 L\'evy process conditioned to stay positive, denoted by $\widehat{X}^{\uparrow}$, is obtained as $h$-transform,
 in the Doob sense, of the killed process with $h(x)=x^{\alpha-1}, x>0$. Thus, it is plain that $\widehat{X}^{\uparrow}$ is the dual, with respect to
the reference measure $y^{\alpha-1}dy$, of $X$. Moreover, the L\'evy
process associated via the Lamperti mapping to
$\widehat{X}^{\uparrow}$ is $(\xi,\P^{(0)})$, see Remark
\ref{rem:l}. Since $X$ and $\widehat{X}^{\uparrow}$ are sspMp, the
conclusion follows from Bertoin and Yor \cite{Bertoin-Yor-02-b}.
\end{proof}

\begin{lemma} \label{lem:ly}
The spectrally positive L\'evy process associated, via the Lamperti
mapping, to the $\kappa$-ssMp $(Y,\Q^{(0)})$ is $(-\xi,\P^{(0)})$.
\end{lemma}
\begin{remark}
As mentioned above, in the case $\kappa=1$, $(-\xi,\P^{(0)})$ corresponds to a Brownian
motion with drift $-\frac{1}{2}$, then $(2Y,\Q^{(0)})$ is a squared
Bessel process of dimension $0$ and we recover the well-known
formula, see e.g.~Revuz and Yor \cite{Revuz-Yor-99},
\begin{eqnarray*}
\Qe^{(0)}_x\left[e^{-2\lambda Y_t}\right]= e^{-x\lambda \left(1+t
2\lambda\right)^{-1}}.
\end{eqnarray*}
\end{remark}
\begin{proof}
First, we show that the branching mechanism associated to
$(Y,\Q^{(0)})$ is the Laplace exponent of a spectrally positive
$\alpha$-stable process.  It is likely that this result is known.
Since we did not find any reference and for sake of completeness we
provide an easy proof. It is well-know, see e.g.~Zeng-Hu \cite{Zeng-06},
that the semi-group of a continuous state branching process with
branching mechanism $\varphi$ admits as spatial Laplace transform
the following expression, for any $\lambda\geq0$,
\begin{eqnarray*}
e^{-x\vartheta_{\lambda}(t)}
\end{eqnarray*}
where $\vartheta:[0,\infty)\rightarrow [0,\infty)$ solves the
following boundary valued differential equation
\begin{eqnarray*}
\vartheta'_{\lambda}(t) =- \varphi(\vartheta_{\lambda}(t)), \quad \vartheta_{\lambda}(0)=\lambda.
\end{eqnarray*}
It is then not difficult to check, in the case
$\varphi(\lambda)=\frac{c}{\kappa}\lambda^{\alpha}$, that
\begin{eqnarray*}
\vartheta_{\lambda}(t) = \lambda \left(1+ct
\lambda^{\kappa}\right)^{-1/\kappa}
\end{eqnarray*}
which is the Laplace exponent of $(Y,\Q^{(0)})$ given in
\eqref{eq:lbs} with $d=1$ and $0<\kappa=\alpha-1\leq 1$. Then, we
deduce that $(Y,\Q^{(0)})$ is obtained from $X$ by the
random time change described above. Finally, from Lemma \ref{lem:fd} and Lemma
\ref{lem:ts}, by choosing $\beta_{\alpha}$ such that
$\beta_{\alpha}=-\frac{\beta_{\alpha}}{\alpha(\beta_{\alpha}-1)}$, i.e.~$\beta_{\alpha}=\frac{\kappa}{\alpha}$,
we get that the image via the Lamperti mapping of $(Y^{\beta_{\alpha}},\Q^{(0)})$ is $(\beta_{\alpha}\xi,\
\P^{(0)})$. By using Lemma
\ref{lem:ps} and observing  that  $(Y,\Q^{(0)})$ is
$\kappa$-self-similar we complete the proof.
\end{proof}
 The proof of the
Theorem \ref{thm:e1} follows readily by observing that for
$\varphi(\lambda)=\frac{c}{\kappa}\lambda^{\alpha}$,
$g(t)=\left(c t\right)^{{\frac{1}{1-\alpha}}}$ and from the
identity
$(i,\Q^{(0)}_x)\stackrel{(d)}{=}(x^{\kappa}\int_0^{\infty}e^{-\kappa
\xi_s}ds, \P_0^{(0)})$. In the following sub-Section, we shall
provide the expression of the semi-group of $(Y,\Q^{(0)})$ in terms of a power
series. We conclude this part by providing the Laplace transform of
the first passage time below for the continuous state branching
process $(Y,\Q^{(0)})$. That is for the stopping time
\begin{eqnarray*}
T^Y_a&=&\inf\{s\geq 0;\: Y_s=a\}.
\end{eqnarray*}
\begin{prop}
For any $x>a>0$, we have
\begin{eqnarray*}
\Qe^{(0)}_x\left[e^{-qT^Y_a}\right]=\frac{\widehat{\mathcal{I}}(q^{\frac{1}{\kappa}}c\kappa
x)}{\widehat{\mathcal{I}}(q^{\frac{1}{\kappa}}c\kappa a)},\quad
q\geq0.
\end{eqnarray*}
where
\begin{eqnarray*}
\widehat{\mathcal{I}}(x) &=& x\int_0^{\infty}e^{-t-xt^{-{\frac{1}{\kappa}}}}t^{-\frac{1}{\kappa}-1}dt \\
&=& \sum_{n=0}^{\infty}
(-1)^{n} \frac{\Gamma(1-\frac{n}{\kappa})}{n!}
x^{n}+\kappa\sum_{n=0}^{\infty}
(-1)^{n} \frac{\Gamma(\kappa(n+1))}{n!}
x^{\kappa(n+1)}\\
&=& \fx{1}{0} \left(
\left.\begin{array}{c}
                  (-\frac{1}{\kappa},1)\\

                \end{array} \right|
 \: -x\right) +\kappa \fx{1}{0} \left(
\left.\begin{array}{c}
                  (\kappa,\kappa)\\

                \end{array} \right|
 \: -x\right)
\end{eqnarray*}
\begin{proof}
First, from \eqref{eq:li}, we deduce readily the identity
\begin{eqnarray*}
\Qe^{(0)}_x\left[e^{-qT^Y_0}\right]=\widehat{\mathcal{I}}(q^{\frac{1}{\kappa}}c\kappa
x),\quad q\geq0.
\end{eqnarray*}
The first part of the claims is completed by  an application of the
strong Markov property and using the fact that $(Y,\Q^{(0)})$ has no
negative jumps. To get the expression of the integral as a power
series we follow a line of reasoning similar to Neretin
\cite{Neretin-06}.  First, consider the space
$L^2(\R^+,\frac{dx}{x})$ and denote by $I_{\varrho,h}(x,v)$ for
$\varrho<0 $ and $\Re(x),\Re(v),\Re(h)>0$ the inner product of the
functions $e^{\varrho}_{x}(t)=e^{-xt^{\varrho}}$ and
$e_{v,h}(t)=t^he^{-vt^{\varrho}}$, i.e.
\begin{eqnarray*}
I_{\varrho,h}(x,v)&=& \int_0^{\infty}e^{-vt-xt^{\varrho}}t^{h-1}dt.
\end{eqnarray*}
Then, the Mellin  transform of $e^{\varrho}_{x}$ is
\begin{eqnarray*}
\tilde{e}^{\varrho}_{x}(\lambda)&=& \int_0^{\infty}e^{-xt^{\varrho}}t^{\lambda-1}dt \\
&=& \frac{\textrm{sgn}(\varrho)}{\varrho}\int_0^{\infty}e^{-xu}u^{\lambda/\varrho-1}du\\
&=& \frac{\textrm{sgn}(\varrho)\Gamma(\lambda/\varrho)}{\varrho
x^{\lambda/\varrho}}.
\end{eqnarray*}
While the Mellin  transform of $e^{v,h}_{x}$ is
\begin{eqnarray*}
\tilde{e}^{v,h}_{x}(\lambda)&=&v^{-h+\lambda}\Gamma(h).
\end{eqnarray*}
By the Plancherel formula for the Mellin transform, we have
\begin{eqnarray*}
\int_{0}^{\infty}e^{\varrho}_{x}\overline{e^{v,h}_{x}}\frac{dx}{x}&=& \frac{1}{2\pi}\int_{-\infty}^{\infty}\tilde{e}^{v,h}_{x}(i\lambda)\overline{\tilde{e}^{v,h}_{x}(i\lambda)}d\lambda
\end{eqnarray*}
that is
\begin{eqnarray*}
I_{\varrho,h}(x,v)&=& \frac{\textrm{sgn}(\varrho)}{2\pi\varrho v^h
}\int_{-\infty}^{\infty}\Gamma(\lambda/\varrho)\Gamma(h+i\lambda)u^{-i\lambda/\varrho}v^{i\lambda}d\lambda.
\end{eqnarray*}
Then we perform the change of variable $z=i\lambda/\varrho$ and consider an arbitrary contour $\mathfrak{C}$  coinciding with the imaginary axis near $\pm \infty$ and leaving all the poles of the integrand of the left side. We get the series expansions by summing the  residues and by choosing $h=v=1$.
\end{proof}
\end{prop}
\subsection{$\ldots$ a family of continuous state branching processes with immigration}

We recall that $\kappa=\alpha-1$, and for any $\delta \in \R^+$ we
write $\P^{(0,\delta)}=\left(\P^{(0,\delta)}_x\right)_{x\in\R}$ for the family of probability measures of the L\'evy process $\xi$, which
admits the following Laplace exponent
\begin{eqnarray*}
\psi^{(0,\delta)}(\lambda)&=&\psi^{(0)}(\lambda)-\frac{\alpha \delta}{\lambda+\kappa}\psi^{(0)}(\lambda),\quad   \lambda \geq 0\\
&=& c\left(\lambda+\kappa-\alpha\delta\right)\frac{\Gamma(\lambda+\kappa)}{\Gamma(\lambda)} .
\end{eqnarray*}
We start by  computing the characteristic
triplet of $(\xi,\P^{(0,\delta)})$.
\begin{prop}
The characteristic triplet of $(\xi,\P^{(0,\delta})$ is given
by   $\sigma^{(\delta)}=0$,
\begin{eqnarray*}
b^{(\delta)}&=& c\Gamma(\alpha)\left(1-\frac{\alpha \delta}{\kappa}\right),
\end{eqnarray*}
and
\begin{eqnarray*}
\nu^{(\delta)}(dy)&=&c_{\alpha}\frac{e^{\alpha
y}}{(1-e^y)^{\alpha+1}}\left(1+\delta(e^{-y}-1)\right) dy, \quad
y<0.
\end{eqnarray*}
where
\begin{eqnarray*}
\log\left(\E^{(0,\delta)}\left[e^{-\lambda
\xi_1}\right]\right)=b^{(\delta)}\lambda
+\int_{-\infty}^0\left(e^{\lambda y}-1-\lambda y\right)
\nu^{(\delta)}(dy).
\end{eqnarray*}
\end{prop}
\begin{proof}
First, writing simply
here $\psi$ for $\psi^{(0)}$, we have from \eqref{eq:if} by choosing
$\gamma=0$ that
\begin{eqnarray*}
\frac{\alpha \delta}{\lambda+\kappa}\psi(\lambda)&=&\frac{c\alpha
\delta}{\lambda+\kappa}\frac{\Gamma(\lambda+\alpha)}{\Gamma(\lambda)}
\\ &=& \frac{c \alpha \delta}{\Gamma(1-\alpha)}\int_{0}^1\left(u^{\lambda}-1\right)u^{\alpha-2}(1-u)^{-\alpha}du\\
&=&-\delta\int_{-\infty}^0\left(e^{\lambda
y}-1\right)\frac{c_{\alpha}e^{\kappa y}dy}{(1-e^{y})^{\alpha}}\\
&=&-\delta\left(\int_{-\infty}^0\left(e^{\lambda y}-1-\lambda
y\right)\frac{c_{\alpha}e^{\kappa
y}dy}{(1-e^{y})^{\alpha}}+\lambda\int_{-\infty}^0
y\frac{c_{\alpha}e^{\kappa y}dy}{(1-e^{y})^{\alpha}}\right)\\
&=&-\delta \int_{-\infty}^0\left(e^{\lambda y}-1-\lambda
y\right)\frac{c_{\alpha}e^{\kappa y}dy}{(1-e^{y})^{\alpha}}+\lambda
c\Gamma(\alpha)\frac{\alpha \delta}{\kappa}.
\end{eqnarray*}
where we have used the identities
$c_{\alpha}=\frac{c}{\Gamma(-\alpha)}
>0$ and
\begin{eqnarray*}
\int_{-\infty}^0 y\frac{c_{\alpha}e^{\kappa
y}dy}{(1-e^{y})^{\alpha}} &=&-\lim_{\lambda \rightarrow 1}
\frac{\partial{}}{\partial \lambda} \int_{-\infty}^0\left(e^{\lambda
y}-1-\lambda y\right)\frac{c_{\alpha}e^{\kappa
y}dy}{(1-e^{y})^{\alpha+1}}\\
&=&-\lim_{\lambda \rightarrow 1} \frac{\partial{}}{\partial \lambda}
\psi^{(-1)}(\lambda)-c\frac{\Gamma(\alpha)}{\kappa}\\
&=& -\psi'(0)-c\frac{\Gamma(\alpha)}{\kappa}.
\end{eqnarray*}
The proof is completed by putting pieces together.
\end{proof}
Next,  set $\delta_{\kappa}=\frac{\delta}{\kappa}$ and
$M_{\delta}=\E^{(0,\delta)}[\xi_1]$ and note from the previous
proposition that
\begin{eqnarray*}
M_{\delta}&=&c\Gamma(\alpha)\left(1-\alpha
\delta_{\kappa}\right).
\end{eqnarray*}
In particular, we have $M_{\delta}<0$ if
$\delta>\frac{\kappa}{\alpha}$. Under such a condition, we write simply
\begin{eqnarray*}
\left(\Sigma_{\infty},\P_0^{(0,\delta)}\right)=\left(\int_0^{\infty}e^{\kappa\xi_s }ds,\P_0^{(0,\delta)}\right).
\end{eqnarray*}
We are now ready to state to the
main result of this section.
\begin{thm} \label{thm:en}
Let $0<\kappa<1$. Then, for any $\delta >\frac{\kappa}{\alpha}$, the law of the positive
random variable $(\Sigma_{\infty}, \P_0^{(0,\delta)})$ is absolutely
continuous with an infinitely continuously differentiable density, denoted by
$f^{(\delta)}_{\infty}$, and
\begin{eqnarray*}
f^{(\delta)}_{\infty}(y)&=& M_{\delta}
y^{-\alpha\delta}\sum_{n=0}^{\infty} (-1)^{n}
\frac{\Gamma(n+\alpha\delta_\kappa)}{n!\Gamma(\kappa
n+\alpha\delta)}
y^{-n}\\
&=&M_{\delta}y^{-\alpha \delta}\fx{1}{1}\left(
\left.\begin{array}{c}
                  (1,\alpha\delta_\kappa)\\
                  (\kappa,\alpha\delta)
                \end{array} \right|
 \: -y^{-1}\right), \: y>0.
\end{eqnarray*}
We deduce, from Section \ref{sec:w},  the following asymptotic behaviors
\begin{eqnarray*}
f^{(\delta)}_{\infty}(y)&\sim& M_{\delta}
\sum_{n=0}^{\infty} (-1)^{n}
\frac{\Gamma(n+\alpha\delta_\kappa)}{n!\Gamma(-\kappa
n)} y^{n} \quad \textrm{ as } y\rightarrow 0.\\
f^{(\delta)}_{\infty}(y)&\sim& M_{\delta}
\frac{\Gamma(\alpha\delta_\kappa)}{\Gamma(\alpha\delta)}y^{-\alpha\delta} \quad \textrm{ as } y\rightarrow \infty.
\end{eqnarray*}

\end{thm}
\begin{remark}
\begin{enumerate}
 \item For $\delta>\frac{\kappa}{\alpha}$, the law of $ \left(\Sigma_{\infty},\P^{(0,\delta)}\right)$ is a generalization of the inverse  Gamma distribution.
Indeed, specifying on $\kappa=1$, i.e.~$\alpha=2$, $c=\frac{1}{2}$
and $\delta>\frac{1}{2}$, the expression above reduces
 to
\begin{eqnarray*}
f^{(\delta)}_{\infty}(y)&=&2
\frac{y^{-2\delta}}{\Gamma(2\delta-1)}e^{-\frac{1}{y}}
\end{eqnarray*}
which corresponds to Dufresne's result \cite{Dufresne-90}, i.e.
\begin{eqnarray*}
\int_0^{\infty}e^{B_s-(\delta-\frac{1}{2})
s}ds&\sim&\frac{1}{2G(2\delta-1)}
\end{eqnarray*}
where we recall that $G(\delta)$ is a Gamma random variable with
parameter $\delta>0$.
\item Note also that, the densities $M^{-1}_{\delta}f_{\infty}^{(\delta)}$ converges, as
$\delta\rightarrow\frac{\kappa}{\alpha}$
 to  a density probability distribution given by
\begin{eqnarray}
f^{(\frac{\kappa}{\alpha})}(y)&=&y^{-1}\sum_{n=0}^{\infty}
 \frac{(-y)^{-n}}{\Gamma(\kappa (n+1))}, \: y>0,
\end{eqnarray}
which is the inverse of a positive Linnik law of parameters
$(\kappa, \kappa)$, see \cite{Linnik-63}.
\end{enumerate}
\end{remark}
The proof of the Theorem above in split in several intermediate
results which we find worth being stated. The plan of the proof is
as follows. We first characterize, in terms of their spatial Laplace
transform the class of continuous state branching processes with
immigration (for short {\emph{cbip}}) which enjoy the scaling property. On
the one hand, by inverting this Laplace transform, we provide an
expression for the density of the semi-group  of this family. In
particular, we shall obtain  the expression of the density of their
entrance laws at $0$. On the other hand, we shall characterize the
L\'evy processes associated to this family of sspMp via the Lamperti
mapping. Finally, we shall show how to relate the density of these
entrance laws to the density of the law of the exponential
functionals under study.

We start with the following easy result which gives a complete
characterization, in terms of their  Laplace transforms, of
self-similar cbip. To this end, we now recall the definition of a
cbip with parameters $[\varphi,\chi]$. It is well known, see
e.g.~\cite{Zeng-06}, that the semi-group of a cbip  with branching
mechanism $\varphi$ and immigration mechanism $\chi$, where $\chi$
is the Laplace exponent of a positive infinitely divisible random
variable, admits
 as a spatial Laplace transform the following expression
\begin{eqnarray} \label{eq:lbis}
\exp\left(-x\vartheta_{\lambda}(t)-\int_0^t\chi(\vartheta_{\lambda}(s))ds\right),\qquad
x,t\geq0.
\end{eqnarray}

\begin{lemma}
A cbip is self-similar of index $\kappa$ if and only if $0<
\kappa\leq1$ and it corresponds to the cbip with parameters
$[\varphi,\delta \chi]$ where $\delta>0$,
$\varphi(\lambda)=\frac{c}{\kappa}\lambda^{\kappa+1}$ and
$\chi(\lambda)=c(\frac{\kappa+1}{\kappa})\lambda^{\kappa}$. Its
Laplace transform has the following expression, for
$\delta,x,\lambda> 0$
\begin{eqnarray*}
\Qe^{(\delta)}_x\left[e^{-\lambda Y_t}\right]= \Lambda^{(\delta)}_{t}(\lambda,x)
\end{eqnarray*}
where
\begin{eqnarray*}
 \Lambda^{(\delta)}_{t}(\lambda,x)&=&\left(1+ct
\lambda^{\kappa}\right)^{-\alpha \delta_{\kappa}}e^{-x\lambda
\left(1+ct \lambda^{\kappa}\right)^{-1/\kappa}}.
\end{eqnarray*}
In particular its entrance law is characterized by
\begin{eqnarray} \label{eq:lte}
\Qe^{(\delta)}_{0^+}\left[e^{-\lambda Y_t}\right]= \left(1+ct
\lambda^{\kappa}\right)^{-\alpha \delta_{\kappa}}.
\end{eqnarray}
We denote this family of processes by
$(Y,\Q^{(\delta)})_{\delta>0}$.
\end{lemma}
\begin{proof}
The sufficient part follows readily from the definition of the cbip
and by observing that for any $a>0$,
$\Lambda^{(\delta)}_{a^{\kappa}t}(\lambda,ax)=\Lambda^{(\delta)}_{t}(a\lambda,x)$.
The necessary part follows from the fact that the unique
self-similar branching process has its Laplace transform given by
\eqref{eq:lbs} and thus the immigration has to satisfy the
self-similarity property. Since we have for any $a>0$,
$a\vartheta_{\lambda}(a^{\kappa}t) =\vartheta_{a\lambda}(t)$, we
need that
\begin{eqnarray*}
\chi(a\vartheta_{a\lambda}(t))=a^{\kappa}\chi(\vartheta_{a\lambda}({\kappa}t))
\end{eqnarray*}
which is only possible for $\chi(a)=Ca^{\kappa}$ for some positive
constant $C$ (since  $\chi$ is the Laplace exponent of a
subordinator). The claim follows.
\end{proof}
\begin{remark}
\begin{enumerate}
\item We mention that $(Y,\Q^{(1)})$ corresponds to the continuous state
branching process $(Y,\Q^{(0)})$ conditioned to never extinct in the terminology of Lambert \cite{Lambert-07}.
Indeed, he showed that the latter corresponds
to the cbpi with immigration $\varphi'$.
\item The Laplace transform of the entrance law \eqref{eq:lte} appears in a paper of Pakes \cite{Pakes-98} where he studies scaled  mixtures of (symmetric) stable laws. More precisely, denoting by $Y_1^{(\delta)}$  the entrance law at time $1$ of $(Y,\Q^{(\delta)})$, we have the following identity in law
\begin{eqnarray*}
Y_1^{(\delta)} \stackrel{(d)}{=} G(\delta)^{\kappa} S_{\kappa},
\end{eqnarray*}
where $S_{\kappa}$ is a positive stable law of index $\kappa$  and the two random variables on the right-hand side are considered to be independent.
\end{enumerate}
\end{remark} We proceed by
providing the expression of the semi-group of
$(Y,\Q^{(\delta)})$.

\begin{prop}
For any $\delta>0$, the semi-group of $(Y,\Q^{(\delta)})$ admits a density,  with respect
to the Lebesgue measure, denoted by $p^{{(\delta)}}_t(.,.)$, which is given for any $x,y,t>0$ by
\begin{eqnarray*}
p^{{(\frac{\kappa\delta}{\alpha})}}_t(x,y)&=&\left(\frac{y}{(ct)^{1/\kappa}}\right)^{\kappa
\delta-1}(ct)^{-1/\kappa}\sum_{n=0}^{\infty}
 \frac{ \fx{1}{1}\left( \left.\begin{array}{c}
                  (1,\frac{n}{\kappa}+\delta) \\
                  (\kappa,\kappa \delta)
                \end{array} \right|
 \:
 -\frac{y^{\kappa}}{ct}\right)}{n!\Gamma\left(\frac{n}{\kappa}+\delta\right)} \left(-\frac{yx}{(ct)^{2/\kappa}}\right)^n.
\end{eqnarray*}
In particular, the density, with respect to the Lebesgue measure, of the self-similar branching process $(Y,\Q^{(0)})$ is given by
\begin{eqnarray*}
p^{(0)}_t(x,y)&=&
\left(\frac{y}{(ct)^{1/\kappa}}\right)^{-1}(ct)^{-1/\kappa}\sum_{n=0}^{\infty}
 \frac{ (-1)^{n}}{n!\Gamma\left(\frac{n}{\kappa}\right)} \fx{1}{1}\left( \left.\begin{array}{c}
                  (1,\frac{n}{\kappa}) \\
                  (\kappa,0)
                \end{array} \right|
 \:
 -\frac{y^{\kappa}}{ct}\right)\left(\frac{yx}{(ct)^{2/\kappa}}\right)^n.
\end{eqnarray*}
For $\delta >0$, the entrance law of $(Y,\Q^{(\delta)})$ is given by
\begin{eqnarray} \label{eq:lbis}
p_t^{^{(\frac{\kappa\delta}{\alpha})}}(0,y)&=&\frac{t^{-1/\kappa}}{\Gamma(\delta)}
\left(\frac{y}{(ct)^{1/\kappa}}\right)^{\kappa\delta-1}
\fx{1}{1}\left( \left.\begin{array}{c}
                  (1, \delta) \\
                  (\kappa,\kappa \delta)
                \end{array} \right|
 \: -\frac{y^{\kappa}}{ct}\right).
\end{eqnarray}
\end{prop}
\begin{proof}
In what follows, we simply write $\delta$ for $\alpha
\delta_{\kappa}$ and set $c=1$. Thus, by means of the binomial
formula, we get, for $t\lambda^{\kappa}>1$,
\begin{eqnarray*}
(1+t\lambda^{\kappa})^{-\delta} &=& \sum_{n=0}^{\infty} (-1)^{n}
\frac{\Gamma(\delta+n)}{n!\Gamma(\delta)}
(t\lambda^{\kappa})^{-n-\delta}.
\end{eqnarray*}
The term-by term inversion yields
\begin{eqnarray*}
p^{(\frac{\kappa\delta}{\alpha})}_{t}(0,y) &=&
\frac{t^{-1/\kappa}}{\Gamma(\delta)}\left(\frac{y}{t^{1/\kappa}}\right)^{\kappa\delta-1}\sum_{n=0}^{\infty}
(-1)^{n} \frac{\Gamma(n+\delta)}{n!\Gamma(\kappa( n+\delta))}
\left(\frac{y}{t^{1/\kappa}}\right)^{\kappa n}.
\end{eqnarray*}
Next, we have
\begin{eqnarray*}e^{-\lambda x(1+t\lambda^{\kappa})^{-\frac{1}{\kappa}}}
(1+t\lambda^{\kappa})^{-\delta}&=& \sum_{n=0}^{\infty}
(-1)^{n}(1+t\lambda^{\kappa})^{-(\frac{n}{\kappa}+\delta)}
\frac{1}{n!} \lambda^nx^n\\
&=& t^{-\delta}\sum_{n=0}^{\infty}
(-1)^{n}(1+(t\lambda^{\kappa})^{-1})^{-(\frac{n}{\kappa}+\delta)}
\frac{1}{n!}
\lambda^{-\kappa\delta}\left(\frac{x}{t^{1/\kappa}}\right)^n.
\end{eqnarray*}
Once again inverting term and term and using the previous result, we deduce that
\begin{eqnarray*}
p^{(\frac{\kappa\delta}{\alpha})}_{t}(x,y)&=&
t^{-\delta}\sum_{n=0}^{\infty} (-1)^{n} F^n_t(y) \frac{1}{n!}
\left(\frac{x}{t^{1/\kappa}}\right)^n
\end{eqnarray*}
where the term $F^n_t(y)$ is given by
\begin{eqnarray*}
F^n_t(y) & =& \frac{ y^{\kappa
\delta-1}}{\Gamma\left(\frac{n}{\kappa}+\delta\right)}
\sum_{r=0}^{\infty}
(-1)^{r}\frac{\Gamma(r+\frac{n}{\kappa}+\delta)}{r!\Gamma(\kappa (r+\delta))}\left(\frac{y}{t^{1/\kappa}}\right)^{\kappa r}\\
& =&\frac{ y^{\kappa
\delta-1}}{\Gamma\left(\frac{n}{\kappa}+\delta\right)}
\fx{1}{1}\left( \left.\begin{array}{c}
                  (1,\frac{n}{\kappa}+\delta) \\
                  (\kappa,\kappa\delta)
                \end{array} \right|
 \: -\frac{y^{\kappa}}{t}\right).
\end{eqnarray*}
Note that $F_t^n(y)=G(n)\star
p^{(\frac{\kappa\delta+n}{\alpha})}_{t}(0,y)$.  Thus, by
putting pieces together
 we get
\begin{eqnarray*}
p^{(\frac{\kappa\delta}{\alpha})}_{t}(x,y)&=&
\left(\frac{y}{t^{1/\kappa}}\right)^{\delta-1}t^{-1/\kappa}\sum_{n=0}^{\infty}
 \frac{ (-1)^{n}}{n!\Gamma\left(\frac{n}{\kappa}+\delta\right)} \fx{1}{1}\left( \left.\begin{array}{c}
                  (1,\frac{n}{\kappa}+\delta) \\
                  (\kappa,\kappa\delta)
                \end{array} \right|
 \: -\frac{y^{\kappa}}{t}\right)\left(\frac{yx}{t^{2/\kappa}}\right)^n
\end{eqnarray*}
which completes the proof.
\end{proof}
Note that in the case $\kappa=1$,  we get
\begin{eqnarray*}
F^n_t(y) & =& \frac{ y^{ \delta-1}}{\Gamma\left(n+\delta\right)}
\sum_{r=0}^{\infty}
(-1)^{r}\frac{\Gamma(r+n+\delta)}{r!\Gamma(r+\delta)}\left(\frac{y}{t}\right)^{r}\\
& =&\frac{ y^{
\delta-1}}{\Gamma\left(n+\delta\right)}
\fx{1}{1}\left( \left.\begin{array}{c}
                  (1,n+\delta) \\
                  (1,\delta)
                \end{array} \right|
 \: -\frac{y^{\kappa}}{t}\right).
\end{eqnarray*}
Thus, by means of the formula for  products of power series, we recover the well-known expression of the density of the semi-group of a Bessel squared process, see e.g.~\cite[p.136]{Borodin-Salminen-02},
\begin{eqnarray*}
p^{(\frac{\delta}{2})}_{\frac{t}{2}}(x,y)&=&\left(\frac{y}{t}\right)^{
\delta-1}t^{-1}e^{-\frac{x+y}{t}}\sum_{n=0}^{\infty}\frac{\left(\frac{xy}{t^2}\right)^n}{n!\Gamma(n+\delta)}\\
&=&\left(\frac{y}{xt}\right)^{\frac{
\delta-1}{2}}t^{-1}e^{-\frac{x+y}{t}}I_{\delta-1}\left(\frac{2\sqrt{xy}}{t}\right),
\end{eqnarray*}
where we recall that
\begin{equation*}
I_{\nu}(x)=\sum_{n=0}^{\infty}\frac{(x/2)^{\nu+2n}}{n!\Gamma(\nu+n+1)}
\end{equation*}
stands for the modified Bessel function of the first kind of index $\nu$, see e.g.~\cite{Lebedev-72}.

We proceed by characterizing the L\'evy processes associated to $(Y,\Q^{(\delta)})$ via the Lamperti mapping.
\begin{prop}
For any $\delta \geq 0$, the L\'evy process associated to the $\kappa$-sspMp $(Y,\Q^{(\delta)})$  is the L\'evy process
 $(-\xi,\P^{(0,\delta)})$.
\end{prop}
\begin{proof}
Let us first consider the case $\delta=1$. Lambert \cite{Lambert-07} showed that $(Y,\Q^{(1)})$
 corresponds to the branching process $(Y,\Q^{(0)})$ conditioned to never extinct,
  which is simply the h-transform in the Doob's sense, with $h(x)=x$, of the minimal process $(Y,\Q^{(0)})$.
   Let us now compute the infinitesimal generator of $(Y,\Q^{(1)})$, denoted by $\mathbf{Q}^{(1)}_+$.
 To this end, let us recall that since the process $X$ does not have negative jumps and has a finite mean,
   its infinitesimal generator, denoted by $ \mathbf{Q}^{\dag}_+$, is given, see also \cite{Caballero-Chaumont-06}, for a smooth function $f$ on $\R^+$ with $f(0)=0$ and any $x>0$, by
\begin{eqnarray*}
 \mathbf{Q}^{\dag}_+ f(x) & = &   \int_{0}^{\infty} \left(f(x+y)-f(x)-yf'(x)\right)\frac{c_{-}}{y^{\alpha +1}}dy\\
 &=&
 x^{-\alpha} \int_{1}^{\infty} \left(f(ux)-f(x)-xf'(x)(u-1)\right)\frac{c_{-}}{(u-1)^{\alpha +1}}du
\end{eqnarray*}
where we have performed the change of variable $y=x(u-1)$. Thus, by
a formula of Volkonskii, see e.g.~Rogers and Williams \cite[III.21]{Rogers-Williams-94}, we
deduce that, for a  function $f$ as above and any $x>0$,
\begin{eqnarray*}
 \mathbf{Q}^{(0)}_+ f(x) & = & x\mathbf{Q}^{\dag}_+ f(x)\\
 & = & x  \int_{0}^{\infty} \left(f(x+y)-f(x)-yf'(x)\right)\frac{c_{-}}{y^{\alpha +1}}dy\\
 &=&
 x^{1-\alpha} \int_{1}^{\infty} \left(f(ux)-f(x)-xf'(x)(u-1)\right)\frac{c_{-}}{(u-1)^{\alpha +1}}du
\end{eqnarray*}
Recalling that for any $x>0$, $\mathbf{Q}^{(0)}_+ h(x)=0$ with
$h(x)=x$ we get, by $h$-transform and for a smooth function $f$ on $\R^+$, that
\begin{eqnarray*}
 \mathbf{Q}^{(1)}_+ f(x) & = & \frac{1}{h(x)}\mathbf{Q}^{(0)}_+ (hf)(x)\\
 & = & x^{1-\alpha} \left(\mathbf{Q}^{(0)}_+ f(x) +\int_{0}^{\infty} \left(f(x+y)-f(x)\right)\frac{c_{-}}{y^{\alpha}}\right)dy\\
 &=& x^{1-\alpha} \left(\mathbf{Q}^{(0)}_+ f(x) +\int_{1}^{\infty} \left(f(ux)-f(x)\right)\frac{c_{-}}{(u-1)^{\alpha}}du\right)\\
 &=& x^{1-\alpha} \left(\mathbf{Q}^{(0)}_+ f(x) +\int_{0}^{1} \left(f\left(\frac{x}{u}\right)-f(x)\right)\frac{c_{-}u^{\alpha-2}}{(1-u)^{\alpha}}du\right)
\end{eqnarray*}
We have already shown, see Lemma \ref{lem:ly},  that the L\'evy
process associated via the Lamperti mapping to $(Y,\Q^{(0)}_+)$ is
$(-\xi,\P^{(0)})$. Next, consider the function
$p_{\lambda}(x)=x^{\lambda}$, with $\lambda<0$ and $x>0$, and note
that $\mathbf{Q}^{1}_+ p_{\lambda}(x)=x^{\lambda-\alpha}
\psi(-\lambda)$. Thus, using the integral \eqref{eq:if} we obtain
that
\begin{eqnarray*}
\int_{0}^{1} \left(u^{-\lambda}-1\right)\frac{c_{\alpha}u^{\alpha-2}}{(1-u)^{\alpha}}du&=&(-\lambda)_{\alpha} \frac{c_{\alpha}\Gamma(1-\alpha)}{-\lambda+\kappa}\\
&=&-c(-\lambda)_{\alpha} \frac{1}{-\lambda+\kappa}
\end{eqnarray*}
Using the recurrence formula of the Gamma function, we deduce that
the image, via the Lamperti mapping, of  $(Y,\Q^{(1)})$ is $(-\xi,\P^{(0,1)})$. The general
case is deduced from the previous one by recalling that for any
$\delta,x,\lambda
>0$, and writing $e_{\lambda}(x)=e^{-\lambda x},\:\lambda \geq0$, we
have, see e.g.~\cite{Zeng-06}, for any $x>0$,
\begin{eqnarray*}
 \mathbf{Q}^{(\delta)}_+ e_{\lambda}(x) & = &
 -e_{\lambda}(x)\left(x\varphi(\lambda)+\delta \chi(\lambda)\right)
\end{eqnarray*}
where $\mathbf{Q}^{(\delta)}_+ $ stands for the infinitesimal generator of $(Y,\Q^{(\delta)})$.
\end{proof}
\begin{remark}
We observe that the process $(-\xi,\P^{(0,1)})$ is equivalent to
$(-\xi,\P^{(-1)})$. This is not really surprising since  as
mentioned in the proof, the process $(Y,\Q^{(1)})$ is obtained
from $(Y,\Q^{(0)})$ by $h$-transform with $h(x)=x$. The
corresponding L\'evy process is thus the $\theta=1$-Esscher
transform of $(-\xi,\P^{(0)})$ which is $(-\xi,\P^{(-1)})$.
\end{remark}
The Theorem is proved by putting pieces together and using the
following easy result.
\begin{lemma}
For any $v>0$ and $\delta>\frac{\kappa}{\alpha}$, we have
\begin{eqnarray*}
f_{\infty}^{(\delta)}(
v)=|\E^{(0,\delta)}[\xi_1]|v^{-\frac{1}{\kappa}}p_1^{(\delta)}(0,v^{-\frac{1}{\kappa}}).
\end{eqnarray*}
\end{lemma}
\begin{proof}
In \cite[Lemma 3.2]{Patie-06c}, the following identity is proved
\begin{eqnarray*}
\E^{(0,\delta)}\left[e^{-qy^{\kappa}\Sigma_{\infty}}\right]=|\E^{(0,\delta)}[\xi_1]|\int_0^{\infty}
e^{-qt}p^{(\delta)}_1(0,yt^{-\frac{1}{\kappa}})y^{1-\kappa}dy.
\end{eqnarray*}
Performing the change of variable $t=y^{\alpha}v$, the proof is
completed by invoking  the injectivity of the Laplace transform.
\end{proof}

Recalling  that for $\delta>\frac{\kappa}{\alpha}$ we have
$\psi(\theta)=0$ where $\theta=\alpha\delta-\kappa$.  We deduce
readily, from Rivero \cite{Rivero-05}, the behavior of
$(Y,\Q^{(\delta)})$ at the boundary point $0$.
\begin{prop}
\begin{enumerate}
\item For $\delta\geq\frac{\kappa}{\alpha}$, $0$ is unattainable.
\item For $\delta<\frac{\kappa}{\alpha}$, $0$ is reached a.s.. Moreover, if $0<\delta<\frac{\kappa}{\alpha}$, the boundary $0$ is
recurrent and reflecting, i.e. there exists an unique recurrent
extension of the minimal process which hits and leaves $0$
continuously a.s. and which is $\kappa$-self-similar on
$[0,\infty)$.
\item For $\delta=0$, the point $0$ is a trap.
\end{enumerate}
\end{prop}
\begin{remark}
Note that in the diffusion case, i.e. $\kappa=1$, there is an
absolute continuity relationship between the family of law
$(\Q^{(\delta)})_{\delta>0}$. More precisely, we have  for any
$\delta,\delta_1>0$
\begin{equation} \label{eq:esscher_s}
d\Q_x^{(\delta_1)}{}_{\mid F_t}=e^{-\delta (Y_t-x)
-\frac{\delta^2}{2} \int_0^{t} Y^{-1}_s ds
}d\Q^{(\delta_1+\delta)}_x{}_{\mid F_t}, \quad x,t>0.
\end{equation}
This follows readily by a time change (via the Lamperti mapping) of
the Cameron-Martin formula which relates the laws of a Brownian
motion with different drifts. In the general case, the Esscher
transform studied earlier does not yield such an absolute continuity
relationship since the family of laws of L\'evy processes
$(\xi,\P^{(0,\delta)})_{\delta\geq0}$ are not related by such a
relationship.

\end{remark}

We proceed by characterizing the class of sspMp which enjoy the
infinite decomposability property introduced by Shiga and Watanabe
\cite{Shiga-Watanabe-73}.  More specifically, let $\mathfrak{D}_M^+$
the Skorohod space of nonnegative valued homogeneous Markov process
with c\`adl\`ag paths. Let $(\Q_x)_{x\geq0}, (\Q^1_x)_{x\geq0}$ and
$(\Q^2_x)_{x>0}$ be three systems of probability measures defined on
$(\mathfrak{D}_M^+,\mathfrak{B}(\mathfrak{D}_M^+))$, where
$\mathfrak{B}(\mathfrak{D}_M^+)$ be the $\sigma$-field on
$\mathfrak{D}_M^+$ generated by the Borel cylinder sets. Then,
define
\begin{eqnarray} \label{eq:sw}
\Q = \Q^1\ast\Q^2
\end{eqnarray}
if and only if for every $x,y\geq0,\:
\Q_{x+y}=\Phi(\Q^1_x\times\Q^2_x)$ where $\Phi$ is the mapping
$\mathfrak{D}_M^+\times \mathfrak{D}_M^+\rightarrow
\mathfrak{D}_M^+$ defined by
\begin{eqnarray} \label{eq:sw}
\Phi(x_1,x_2) = x_1+x_2.
\end{eqnarray}
A positive valued Markov process with law $\Q$ is infinitely
decomposable if for every $n\geq 1$ there exits a $\Q^{n}$ such
that
\begin{eqnarray} \label{eq:sw}
\Q = \underbrace{\Q^{n}\ast \ldots \ast \Q^{n}}_{n }.
\end{eqnarray}
They showed that there is one to  one mapping between cbip and
conservative  Markov processes having the property \eqref{eq:sw}.
Thus, we get the following result.
\begin{cor}
There is one to one mapping between the family of sspMp satisfying
the Shiga-Watanabe infinite decomposability property \eqref{eq:sw}
and the family $(Y,\Q^{(\delta)})_{\delta> 0}$, i.e. the family of
spectrally positive sspMp associated, via the Lamperti mapping, to
the family of L\'evy processes $(-\xi,\P^{(0,\delta)})_{\delta> 0}$.
\end{cor}

Finally, we provide an extension of the previous result the
$\kappa$-sspMp Ornstein-Uhlenbeck processes. More specifically,  for
any $\eta\in\R$, let us introduce the family of laws
$(\Q^{\eta,(\delta)})_{\delta\geq0}$ of self-similar
Ornstein-Uhlenbeck processes associated to
$(Y,\Q^{(\delta)})_{\delta\geq0}$ by the following time-space transform, for
any $t\geq0$ and $\delta\geq0$,
\begin{eqnarray*}
\left(Y_t,\Q^{\eta,(\delta)}\right)=\left(e^{-\eta
t}Y_{\tau_{-\eta}(t)},\Q^{(\delta)}\right),
\end{eqnarray*}
where
\begin{eqnarray*}
\tau_{\eta}(t) =\frac{1-e^{-\eta \kappa t}}{\eta \kappa}.
\end{eqnarray*}
For any $\eta \in \R$, $(Y,\Q^{\eta,(\delta)})$ is an homogenous Markov
process and  for $\eta>0$, it has an unique stationary measure which
is the entrance law of $(Y,\Q^{(\delta)})$, see
e.g.~\cite{Patie-ouq-06}. Moreover, its semi-group is characterized by
its Laplace transform as follows
\begin{eqnarray*}
\Qe_x\left[e^{-\lambda U_t}\right]= \left(1+c\tau_{\eta}(t)
\lambda^{\kappa}\right)^{-\delta/\kappa}e^{-xe^{-\eta t}\lambda
\left(1+c\tau_{\eta}(t) \lambda^{\kappa}\right)^{-1/\kappa}},\
\end{eqnarray*}
It is easily  shown that the infinitesimal generator of
$(Y,\Q^{\eta,(\delta)})$ has the following form
\begin{eqnarray*}
 \mathbf{Q}^{\eta,(\delta)}_+ f(x)=  \mathbf{Q}^{(\delta)}_+ f(x) -\eta x f'(x)
\end{eqnarray*}
for a smooth function $f$. Hence we deduce from the identity
\begin{eqnarray*}
 \mathbf{Q}^{(\delta)}_+ e_{\lambda}(x) & = &
 -e_{\lambda}(x)\left(x\varphi(\lambda)+\delta \chi(\lambda)\right),
\end{eqnarray*}
 that $(Y,\Q^{\eta,(\delta)})$ is a cbip with branching
mechanism $\varphi_{\eta}(\lambda)=\varphi{\lambda}-\eta \lambda$
and the immigration mechanism $\chi(\lambda)$. Its semi-group is
absolutely continuous with a density denoted
$p^{{(\delta,\eta)}}_t(x,y)$ and given, for any $x,y,t>0$, by
\begin{eqnarray*}
p^{{(\delta,\eta)}}_t(x,y) &=& e^{-\kappa\eta t}
p^{{(\delta)}}_{\tau_{-\eta}(t)}(x,e^{-\kappa\eta t}y).
\end{eqnarray*}

\section{Some concluding remarks}

\subsection{Representations of some $\fx{p}{q}$  functions}
Let us recall that, for
$\frac{\kappa}{\alpha}<\delta<\frac{2\alpha-2}{\alpha}$, the Laplace
transform of $(\Sigma_{\infty},\P^{(0,\delta)})$ has been computed
by Patie \cite{Patie-06c} as follows, for any $x\geq0$,
\begin{eqnarray} \label{eq:le}
\E^{(0,\delta)} \left[e^{-x \Sigma_{\infty}} \right]
&=&\N_{\kappa,\delta}(x).
\end{eqnarray}
where, by setting $0<m_{\kappa}=2-\frac{\alpha \delta}{\kappa}<1$,
\begin{equation*}
\N_{\kappa,\delta}(x) =\Ip\left(x\right)-
C_{m_{\kappa}}x^{\frac{\alpha
\delta}{\kappa}-1}\mathcal{I}_{\kappa,\delta,m_{\kappa}}\left(x\right),\quad
x\geq0,
\end{equation*}
and
\begin{eqnarray*}
\Ip\left(ckx\right)&=&\Gamma(m_{\kappa})\Gamma(\kappa)
\sum_{n=0}^{\infty}\frac{x^n}{\Gamma(n+
m_{\kappa})\Gamma(\kappa(n+1))}\\
&=& \Gamma(m_{\kappa})\Gamma(\kappa) \fx{1}{2}\left(
\left.\begin{array}{c} (1,1)\\
                  (1,m_{\kappa})(\kappa,\kappa)
                \end{array} \right|
 \: x\right),\\
\mathcal{I}_{\kappa,\delta,\theta}\left(ckx\right)&=&\Gamma(\alpha\delta)
\sum_{n=0}^{\infty}\frac{x^n}{n!\Gamma(\kappa n+\alpha\delta)}\\
&=&\Gamma(\alpha\delta) \quad \fx{0}{1}\left(
\left.\begin{array}{c} (\kappa,\alpha\delta)
                \end{array} \right|
 \: x\right),
\end{eqnarray*}
and where $C_{m_{\kappa}}$ is determined by
\begin{equation} \label{eq:c}
\Ip\left(x\right) \sim
C_{m_{\kappa}}x^{\frac{\alpha
\delta}{\kappa}-1}\mathcal{I}_{\kappa,\delta,m_{\kappa}}\left(x\right)\quad
\textrm{ as } x\rightarrow \infty.
\end{equation}
Using the exponentially infinite asymptotic expansions \eqref{eq:a12} and \eqref{eq:a01}, we deduce that
\begin{equation*}
C_{m_{\kappa}}=\frac{\Gamma(m_{\kappa})\Gamma(\kappa)}{\Gamma(\alpha\delta)}.
\end{equation*}
Observing that $\lim_{\lambda \rightarrow \infty}
\frac{\psi(\lambda)}{\lambda^{\kappa+1}}=c$, we also have, see
\cite{Patie-06c},
\begin{eqnarray*}
 C_{m_{\kappa}} =  \frac{\Gamma(m_{\kappa})}{\kappa
}c^{(m_{\kappa}-1)}e^{E_{\gamma} \kappa
(m_{\kappa}-1)}\prod_{k=1}^{\infty}e^{-\frac{\kappa+\alpha\delta}{k}}\frac{(k+1-m_{\kappa})\psi^{(0,\delta)}(\kappa
k)}{k\psi^{(0,\delta)}(\kappa k+1-m_{\kappa})}
\end{eqnarray*}
where $E_{\gamma}=0.577\ldots$ stands for Euler-Mascheroni constant.

As a consequence of the Theorem \ref{thm:en} we have this
interesting representation of the function $\N_{\kappa,\delta}(x)$
under consideration.
\begin{cor}
For any $\delta
>\frac{\kappa}{\alpha}$, the density $f^{(\delta)}_{\infty}$ is a mixture of exponential
distribution, in particular is log-convex. Moreover, for
$\frac{\kappa}{\alpha}<\delta<\frac{2\alpha-2}{\alpha}$, we have the
following two representations,
\begin{eqnarray*}
\N_{\kappa,\delta}(x) &=&
\frac{M_{\delta}}{\Gamma(\alpha\delta_{\kappa})}
\int_0^{\infty}\frac{u^{\alpha\delta_{\kappa}-1}}{x+u}\quad \fx{0}{1}\left(
\left.\begin{array}{c}
                  (\kappa,\alpha\delta)
                \end{array} \right|
 \: -u\right)du. \\\N_{\kappa,\delta}(x)  &=&
\exp\left(-\int_0^{\infty} (1-e^{-ux})L(u)du\right)
\end{eqnarray*}
where  $L(u)=\int_0^{\infty}e^{-ur}q(r)dr$ with $0\leq q(r)\leq1$ a
measurable function and $\int_0^1\frac{q(r)}{r}dr<\infty$. Finally, we have the following asymptotic expansion for large $x$,
\begin{equation*}
\N_{\kappa,\delta}(x) \sim  \exp\left(-(\kappa^{-\kappa}\alpha^{\alpha}x)^{\frac{1}{\alpha}}\right)x^{-\frac{\delta}{\kappa}-\frac{1}{\alpha}+1}.
\end{equation*}

\end{cor}
\begin{proof}
Note that, for any $\delta >\frac{\kappa}{\alpha}$ and $y>0$, we
have
\begin{eqnarray*}
f^{(\delta)}_{\infty}
(y)&=&\frac{M_{\delta}}{\Gamma(\alpha\delta_{\kappa})}
 \int_0^{\infty}e^{-uy}\sum_{n=0}^{\infty}
\frac{(-1)^{n}}{n!\Gamma(\kappa n+\alpha\delta)}
u^{n+\alpha\delta_{\kappa}-1}du\\
&=&\frac{M_{\delta}}{\Gamma(\alpha\delta_{\kappa})}
\int_0^{\infty}e^{-uy}u^{\alpha\delta_{\kappa}-1}\quad \fx{0}{1}\left(
\left.\begin{array}{c}
                  (\kappa,\alpha\delta)
                \end{array} \right|
 \: -u\right)du.
\end{eqnarray*}
Thus $f^{(\delta)}_{\infty}$ is a completely monotone function. The
fact that $f^{(\delta)}_{\infty}$ is a mixture of exponential
distribution follows from Proposition 51.8  in Sato \cite{Sato-99}.
The first representation is obtained by taking the Laplace transform
on both sides on the previous equation and using \eqref{eq:le}. The
last one follows from \cite[Theorem 51.12]{Sato-99}. The asymptotic behavior of the function $\N_{\kappa,\delta}$ is deduced from the exponentially infinite asymptotic expansions \eqref{eq:a12} and \eqref{eq:a01}.
\end{proof}
Next, since the law of $
\left(\Sigma_{\infty},\P^{(0,\delta)}\right)$ is self-decomposable,
see e.g.~\cite{Patie-06c}, it is unimodal. In other words, the
derivative of the continuously differentiable density admits a unique zero.
 We state this fact in the following.
\begin{cor}
For any $\delta>\frac{\kappa}{\alpha}$ with $0<\kappa\leq 1$, the
function
\begin{eqnarray*}
x\mapsto \fx{1}{1}\left( \left.\begin{array}{c}
                  (1,1+\alpha\delta_\kappa)\\
                  (\kappa,\alpha\delta)
                \end{array} \right|
 \: -x\right)
\end{eqnarray*}
admits a unique zero  on $\R^+$.
\end{cor}
\begin{remark}
In the case $\kappa=1$, the result above is obvious since
\begin{eqnarray*}
 \fx{1}{1}\left( \left.\begin{array}{c}
                  (1,1+\alpha\delta)\\
                  (1,\alpha\delta)
                \end{array} \right|
 \: -x\right) = e^{-x}\left(\alpha\delta-x\right).
\end{eqnarray*}

\end{remark}

\section{Asymptotic expansions of the Wright Hypergeometric
functions} \label{sec:w} Special classes of the Wright generalized
hypergeometric functions have been considered among others by
Mittag-Leffler \cite{Mittag-03}, Barnes \cite{Barnes-07}, Fox
\cite{Fox-28} while the general case $\fx{p}{q}$ has been considered
by  Wright \cite{Wright-40}. We refer to Braaksma \cite[Chap.
12]{Braaksma-64} for a detailed account of this function and its
relation to the $G$-function. In the sequel, we simply indicate
special properties, which can be found in \cite[Chap.
12]{Braaksma-64}.

We proceed by recalling that the Wright
hypergeometric function is defined as
\begin{equation*}
 \quad \fx{p}{q}\left( \left.\begin{array}{c}
                  (A_1,a_1),\ldots(A_p,a_p) \\
                  (B_1,b
                  _1)\ldots(B_q,b_q)
                \end{array} \right|
 \: z\right)
= \sum_{n=0}^{\infty} \frac{\prod_{i=1}^p\Gamma(A_i n+a_i
)}{\prod_{i=1}^q\Gamma(B_in +b_i)}\frac{z^n}{n!}
\end{equation*}
where $p,q$ are nonnegative integers, $a_i \in \C \: (i=1\ldots
p), b_j\in \C \: (j=1\ldots q)$, the coefficients $A_i \in
\R^+\: (i=1\ldots p)$ and $B_j \in \R^+ \:(j=1\ldots q)$ are such
that
\begin{equation*}
A_i n+a_i \neq 0,-1,-2,\ldots \qquad(i=0,1,\ldots,p;n=0,1,\ldots).
\end{equation*}
In what follows, we will also use the number $S$ and $T$ defined respectively by
\begin{eqnarray*}
S &=& 1+\sum_{i=1}^q B_i-\sum_{i=1}^p A_i\\
T &=& \prod_{i=1}^pA_i^{A_i}\prod_{i=1}^q B_i^{-B_i}.
\end{eqnarray*}
Throughout this part we assume that $S$ is positive. In such a case,
the series is convergent for all values of $z$ and it defines an
integral function of $z$ (the case  $S=0$ is also treated in
\cite{Braaksma-64}). Next, for $S$ positive,  the function  $
\fx{p}{q}$ admits a contour integral representation. More precisely,
we have

\begin{equation*}
 \quad \fx{p}{q}\left( \left.\begin{array}{c}
                  (A_1,a_1),\ldots(A_p,a_p) \\
                  (B_1,b
                  _1)\ldots(B_q,b_q)
                \end{array} \right|
 \: z\right)
= \frac{1}{2\pi i} \int_{\mathfrak{C}} \frac{\prod_{i=1}^p\Gamma(A_i n+a_i
)}{\prod_{i=1}^q\Gamma(B_in +b_i)}\Gamma(-s)(-s)^zds
\end{equation*}
where $\mathfrak{C}$ is a contour in the complex $s$-plane which
runs from $s=a-i\infty$ to $s=a+i\infty$ ($a$ an arbitrary real
number) so that the points $s=0,1,2 \ldots$ and
$s=-\frac{a_j+n}{A_i}, \: (i=0,1,\ldots,p;n=0,1,\ldots) $ lie to the
right left of $\mathfrak{C}$. Next we introduce the following
functions
\begin{eqnarray*}
P(z)&=& \sum_{s\in R_p} z^s \Gamma(-s) Res\left(\frac{\prod_{i=1}^p\Gamma(A_i s+a_i
)}{\prod_{i=1}^q\Gamma(B_i s +b_i)}\right) \\
E(z)&=&
\frac{\exp\left((TS^Sz)^{\frac{1}{S}}\right)}{S}\sum_{k=0}^{\infty}
H_k(TS^Sz)^{\frac{1-G-k}{S}}
\end{eqnarray*}
where $Res$ stands for residuum, we set $R_p=\{r_{i,n}=-\frac{a_i+n}{A_i},i=0,1,\ldots,p;n=0,1,\ldots\}$ and the constant $G$ is given by
\begin{eqnarray*}
G &=& \sum_{i=1}^q b_i-\sum_{i=1}^p a_i+\frac{p-q}{2}+1.
\end{eqnarray*}
The coefficients $ (H_k)_{k\geq0}$ are determined by
\begin{eqnarray*}
\frac{\prod_{i=1}^p\Gamma(A_i s+a_i
)}{\prod_{i=1}^q\Gamma(B_i s +b_i)}(TS^S)^{-s}\sim \sum_{k=0}^{\infty} \frac{H_k}{\Gamma(k+Ss+G)}.
\end{eqnarray*}
In particular,
\begin{eqnarray*}
H_0=(2\pi)^{\frac{p-q}{2}}S^{G-\frac{1}{2}}\prod_{i=1}^pA_i^{a_i-\frac{1}{2}}\prod_{i=1}^qB_i^{\frac{1}{2}-b_i}.
\end{eqnarray*}

We have the following asymptotic expansions.
\begin{enumerate}
\item Suppose $S>0$ and  $p>0$. Then, the following algebraic asymptotic expansion
\begin{eqnarray*}
\quad \fx{p}{q}(z) \sim P(-z)
\end{eqnarray*}
holds for $\mid z\mid \rightarrow \infty$ uniformly on every closed subsector of
\begin{eqnarray*}
\mid \textrm{arg}(-z) \mid< \left(1-\frac{S}{2}\right)\pi
\end{eqnarray*}
\item Suppose $S>0$. Then, the following exponentially infinite asymptotic expansion
\begin{eqnarray*}
\quad \fx{p}{q}(z) \sim E(z)
\end{eqnarray*}
holds for $\mid z\mid \rightarrow \infty$ uniformly on every closed sector (vertex in $0$) contained in $\textrm{arg}(z) < \min(S,2)\frac{\pi}{2}$.
\end{enumerate}
For the convenience of the reader, we list below the asymptotic expansion corresponding to $ \fx{p}{q}$
 functions which  appear in this paper. For $ y \rightarrow \infty$, we have
\begin{eqnarray}
 \fx{1}{1}\left( \left.\begin{array}{c}
                  (1, \alpha\delta_{\kappa}) \\
                  (\kappa,\alpha \delta)
                \end{array} \right|
 \: -y^{\kappa}\right) &\sim & \sum_{n=0}^{\infty} (-1)^n\frac{\Gamma(n+\alpha\delta_{\kappa})
}{\Gamma(-\kappa n)}\frac{y^{-\alpha\delta_{\kappa}-n}}{n!}, \label{eq:a11}\\
 \fx{1}{2}\left(
\left.\begin{array}{c} (1,1)\\
                  (1,m_{\kappa})(\kappa,\kappa)
                \end{array} \right|
 \: y\right) &\sim& (2\pi\alpha)^{-\frac{1}{2}}
 \kappa^{\frac{1}{2}-\delta-\frac{\kappa}{2 \alpha}}
 \exp\left((\kappa^{-\kappa}\alpha^{\alpha}y)^{\frac{1}{\alpha}}\right)y^{\frac{\delta}{\kappa}+\frac{1}{\alpha}-1} \label{eq:a12} \\
 \fx{0}{1}\left(
\left.\begin{array}{c} (\kappa,\alpha\delta)
                \end{array} \right|
 \: y\right) &\sim& (2\pi\alpha)^{-\frac{1}{2}} \kappa^{\frac{1}{2}-\delta -\frac{\kappa}{2\alpha}}
  \exp\left((\kappa^{-\kappa}\alpha^{\alpha}y)^{\frac{1}{\alpha}}\right)y^{\frac{1}{2\alpha}-\delta}. \label{eq:a01}
\end{eqnarray}


\begin{thebibliography}{10}

\bibitem{Barnes-07}
E.W. Barnes.
\newblock The asymptotic expansion of integral functions defined by generalized
  hypergeometric series.
\newblock {\em Proc. London Math. Soc.}, 5(2):59--116, 1907.

\bibitem{Bertoin-96}
J.~Bertoin.
\newblock {\em L\'evy Processes}.
\newblock Cambridge University Press, Cambridge, 1996.

\bibitem{Bertoin-Biane-Yor-04}
J.~Bertoin, P.~Biane, and M.~Yor.
\newblock Poissonian exponential functionals, {$q$}-series, {$q$}-integrals,
  and the moment problem for log-normal distributions.
\newblock In {\em Seminar on Stochastic Analysis, Random Fields and
  Applications {IV}}, volume~58 of {\em Progr. Probab.}, pages 45--56.
  Birkh\"auser, Basel, 2004.

\bibitem{Bertoin-Yor-02-b}
J.~Bertoin and M.~Yor.
\newblock The entrance laws of self-similar {M}arkov processes and exponential
  functionals of {L}\'evy processes.
\newblock {\em Potential Anal.}, 17(4):389--400, 2002.

\bibitem{Bertoin-Yor-05}
J.~Bertoin and M.~Yor.
\newblock {Exponential functionals of L\'evy processes}.
\newblock {\em Probab. Surv.}, 2:191--212, 2005.

\bibitem{Borodin-Salminen-02}
A.N. Borodin and P.~Salminen.
\newblock {\em Handbook of Brownian Motion - Facts and Formulae}.
\newblock Probability and its Applications. Birkh\"auser Verlag, Basel,
  $2^{nd}$ edition, 2002.

\bibitem{Braaksma-64}
B.L.J. Braaksma.
\newblock Asymptotic expansions and analytic continuations for a class of
  {B}arnes-integrals.
\newblock {\em Compositio Math.}, 15:239--341, 1964.

\bibitem{Caballero-Chaumont-06}
M.E. Caballero and L.~Chaumont.
\newblock Conditioned stable {L}\'evy processes and the {L}amperti
  representation.
\newblock {\em J. Appl. Probab.}, 43(4):967--983, 2006.

\bibitem{Carmona-Petit-Yor-01}
Ph. Carmona, F.~Petit, and M.~Yor.
\newblock Exponential functionals of {L}\'evy processes.
\newblock {\em L\'evy processes Theory and Applications}, pages 41--55, 2001.

\bibitem{Duffie-Filipovic-03}
D.~Duffie, D.~Filipovi\'c, and W.~Schachermayer.
\newblock Affine processes and applications in finance.
\newblock {\em Ann. Appl. Probab.}, 13(3):984--1053, 2003.

\bibitem{Dufresne-90}
D.~Dufresne.
\newblock The distribution of a perpetuity, with applications to risk theory
  and pension funding.
\newblock {\em Scand. Actuar. J.}, (1-2):39--79, 1990.

\bibitem{Erdelyi-55}
A.~Erd{\'e}lyi, W.~Magnus, F.~Oberhettinger, and F.G. Tricomi.
\newblock {\em Higher Transcendental Functions}, volume~3.
\newblock McGraw-Hill, New York-Toronto-London, 1955.

\bibitem{Fox-28}
C.~Fox.
\newblock The asymptotic expansion of generalized hypergeometric functions.
\newblock {\em Proc. London Math. Soc.}, 27(2):389--400, 1928.

\bibitem{Getoor-Sharpe-84}
R.K. Getoor and M.J. Sharpe.
\newblock Naturality, standardness, and weak duality for {M}arkov processes.
\newblock {\em Z. Wahrsch. verw. Gebiete}, 67:1--62, 1984.

\bibitem{Gradshteyn-Ryzhik-00}
I.S. Gradshteyn and I.M. Ryshik.
\newblock {\em Table of Integrals, Series and Products}.
\newblock Academic Press, San Diego, $6^{th}$ edition, 2000.

\bibitem{Grey-74}
D.R. Grey.
\newblock Asymptotic behaviour of continuous time, continuous state-space
  branching processes.
\newblock {\em J. Appl. Probability}, 11:669--677, 1974.

\bibitem{Lambert-07}
A.~Lambert.
\newblock Quasi-stationary distributions and the continuous-state branching
  process conditioned to be never extinct.
\newblock {\em Electron. J. Probab.}, 12:420--446, 2007.

\bibitem{Lamperti-67a}
J.~Lamperti.
\newblock Continuous state branching processes.
\newblock {\em Bull. Amer. Math. Soc.}, 73:382--386, 1967.

\bibitem{Lamperti-67b}
J.~Lamperti.
\newblock Limiting distributions for branching processes.
\newblock In {\em Proc. Fifth Berkeley Sympos. Math. Statist. and Probability
  (Berkeley, Calif., 1965/66), Vol. {II}: Contributions to Probability Theory,
  Part}, pages 225--241. Univ. California Press, Berkeley, Calif., 1967.

\bibitem{Lamperti-72}
J.~Lamperti.
\newblock Semi-stable {M}arkov processes. {I}.
\newblock {\em Z. Wahrsch. Verw. Geb.}, 22:205--225, 1972.

\bibitem{Lebedev-72}
N.N. Lebedev.
\newblock {\em Special Functions and their Applications}.
\newblock Dover Publications, New York, 1972.

\bibitem{Linnik-63}
Ju.~V. Linnik.
\newblock Linear forms and statistical criteria. {I}.
\newblock In {\em Selected Transl. Math. Statist. and Prob., Vol.}, pages
  1--40. Amer. Math. Soc., Providence, R.I., 1963.

\bibitem{Matsumoto-Yor-05-1}
H.~Matsumoto and M.~Yor.
\newblock Exponential functionals of {B}rownian motion. {I}. {P}robability laws
  at fixed time.
\newblock {\em Probab. Surv.}, 2:312--347, 2005.

\bibitem{Matsumoto-Yor-05-2}
H.~Matsumoto and M.~Yor.
\newblock Exponential functionals of {B}rownian motion. {II}. {S}ome related
  diffusion processes.
\newblock {\em Probab. Surv.}, 2:348--384, 2005.

\bibitem{Milgram-06}
M.S. Milgram.
\newblock {On Hypergeometric 3F2(1)}.
\newblock {\em available as http://www.arXiv.org:math:CA/0603096}.

\bibitem{Mittag-03}
G.~Mittag-Leffler.
\newblock Sur la nouvelle function ${E}_{\alpha}(x)$.
\newblock {\em C. R. Math. Acad. Sci. Paris}, 137:554--558, 1903.

\bibitem{Neretin-06}
Y.A. Neretin.
\newblock Stable densities and operators of fractional differentiation.
\newblock In {\em Representation Theory, Dynamical Systems, and Asymptotic
  Combinatorics}, volume 217 of {\em Amer. Math. Soc. Transl. Ser.}, pages
  117--137. Amer. Math. Soc., Providence, RI, 2006.

\bibitem{Pakes-98}
A.G. Pakes.
\newblock Mixture representations for symmetric generalized {L}innik laws.
\newblock {\em Statist. Probab. Lett.}, 37(3):213--221, 1998.

\bibitem{Patie-06c}
P.~Patie.
\newblock Infinitely divisibility of solutions to some semi-stable
  integro-differential equations and exponential functionals of {L\'evy}
  processes.
\newblock {\em Submitted}, 2007.

\bibitem{Patie-ouq-06}
P.~Patie.
\newblock $q$-invariant functions associated to some generalizations of the
  {Ornstein-Uhlenbeck} semigroup.
\newblock {\em Under revision for ALEA}, 2007.

\bibitem{Revuz-Yor-99}
D.~Revuz and M.~Yor.
\newblock {\em Continuous Martingales and Brownian Motion}, volume 293.
\newblock Springer-Verlag, Berlin-Heidelberg, $3^{rd}$ edition, 1999.

\bibitem{Rivero-05}
V.~Rivero.
\newblock Recurrent extensions of self-similar {M}arkov processes and
  {C}ram\'er's condition.
\newblock {\em Bernoulli}, 11(3):471--509, 2005.

\bibitem{Rogers-Williams-94}
L.C.G. Rogers and D.~Williams.
\newblock {\em Diffusions, {M}arkov Processes, and Martingales. {V}ol. 1.}
\newblock Cambridge Mathematical Library. Cambridge University Press,
  Cambridge, 2000.
\newblock Foundations, Reprint of the second (1994) edition.

\bibitem{Sam-Taq-94}
G.~Samorodnitsky and M.S. Taqqu.
\newblock {\em Stable Non-{G}aussian Random Processes}.
\newblock Stochastic Modeling. Chapman \& Hall, New York, 1994.

\bibitem{Sato-99}
K.~Sato.
\newblock {\em L\'evy Processes and Infinitely Divisible Distributions}.
\newblock Cambridge University Press, Cambridge, 1999.

\bibitem{Shiga-Watanabe-73}
T.~Shiga and S.~Watanabe.
\newblock Bessel diffusions as a one-parameter family of diffusion processes.
\newblock {\em Z. Wahrscheinlichkeitstheorie und Verw. Gebiete}, 27:37--46,
  1973.

\bibitem{Wright-40}
E.M. Wright.
\newblock The asymptotic expansion of the generalized hypergeometric function.
\newblock {\em Proc. London Math. Soc. (2)}, 46:389--408, 1940.

\bibitem{Zeng-06}
L.~Zeng-Hu.
\newblock Branching processes with immigration and related topics.
\newblock {\em Front. Math. China}, 1(1):73--97, 2006.

\end{thebibliography}
\end{document}